\documentclass[10pt]{amsart}
\usepackage{files/ehastyle}

\begin{document}
\title{Eckmann-Hilton arguments in equivariant higher algebra}
\author{Natalie Stewart}
\date{\today}

\begin{abstract}
  Let $\cO^{\otimes}$ and $\cP^{\otimes}$ be $k$- and $\ell$-connected unital $G$-operads subject to the condition for all $S$ that $\cO(S) = \emptyset$ if and only if $\cP(S) = \emptyset$.
  We show that the Boardman-Vogt tensor product $\cO^{\otimes} \obv \cP^{\otimes}$ is $(k + \ell + 2)$-connected;
  equivalently, $\cO \otimes \cP$-monoids in any $(k + \ell + 3)$-category lift uniquely to incomplete semi-Mackey functors.
  As a consequence, we show that the smashing localizations on unital $G$-operads correspond precisely to unital $\cN_\infty$-operads, and hence to the (finite) poset of unital weak indexing systems by previous work of the author.
  Along the way we characterize $\ell$-connectivity of a unital $G$-operad $\cO^{\otimes}$ equivalently as $\ell$-connectivity of $\cO$-admissible Wirthm\"uller maps of $\cO$-monoid spaces.
  
  In the discrete case, under no connectivity assumptions, $\cO \otimes \cP$-monoids lift uniquely to incomplete semi-Mackey functors, recovering an Eckmann-Hilton argument for ``$C_p$-unital magmas.'' 
  In the limiting case of infinite tensor powers, we \emph{take the loops out of equivariant infinite loop space theory,} constructing algebraic approximations to incompletely stable $G$-spectra over arbitrary transfer systems.
\end{abstract}

\maketitle

\toc

\newpage
\section*{Introduction}\stoptocwriting
The classical \emph{Eckmann-Hilton argument} shows that, given a set $M$ with two multiplications $*,\cdot\colon M^2 \rightarrow M$ possessing a shared multiplicative unit and 
satisfying the interchange law 
\[
  (a * b) \cdot (c*d) = (a \cdot c) * (b \cdot d),
\]
the unital magmas $(M,*)$ and $(M,\cdot)$ are isomorphic to each other and are commutative monoids.
Indispensable to algebraic topologists, this fact recovers the usual proof that $\pi_n(X)$ is abelian for all pointed spaces $X$ and numbers $n \geq 2$, as well as the same claim for $n = 1$ when $X$ admits an $H$-space structure.
We will study equivariant variations of this result, beginning with a weakening of Dress' Mackey functors \cite{Dress}.
\begin{definition}
  Let $\cC$ be a 1-category with finite products and $C_p$ the cyclic group of prime order $p$.
  A \emph{$C_p$-unital magma in $\cC$} is a unital magma $M^{e}$ with a $C_p$ action by unital magma homomorphisms, a unital magma $M^{C_p}$ (with trivial $C_p$-action), 
  and $C_p$-equivariant \emph{restriction and transfer} homomorphisms
  \[
    r\colon M^{C_p} \rightarrow M^e, \hspace{50pt} t\colon M^e \rightarrow M^{C_p}
  \]
  subject to the condition that $r \circ t$ is multiplication by $p$.
  A \emph{homomorphism} $M \rightarrow N$ is a pair of unital magma homomorphisms $F^e\colon M^e \rightarrow N^e$ and $F^{C_p}\colon M^{C_p} \rightarrow M^e$ such that $F^{C_p} \circ t = t \circ F^e$ and $F^e \circ r = r \circ F^{C_p}$.
\end{definition}
\begin{example}\label{Homotopy coefficient system}
  Let $\rho$ be the regular representation of $C_p$.
  The $\rho$th homotopy coefficient system of a pointed $C_p$-space attains a natural $C_p$-unital magma structure under the evident generalization of Lewis' \emph{$\rho$-Mackey structure} \cite{Lewis_hurewicz}.\footnote{Explicitly, by $V$-Mackey functor, we mean a functor $\sB_G(V) \rightarrow \Ab$ sending disjoint unions to direct sums, where $\sB_G(V)$ is the category whose objects are finite $G$-sets and whose hom sets are $[\Sigma^V T_+, \Sigma^V S_+]$;
  the transfer map $\Sigma^{\rho}_+ *_{C_p} \rightarrow \Sigma^{\rho}_+ [C_p/e]$ is constructed by the usual $\mathbb{S}_G$-duality construction along an embedding $[C_p/e] \hookrightarrow \rho$ (see \cite{Wirthmuller}).}
\end{example}
In this article, we vastly generalize the following theorem, which is proved as \cref{GEHA-precise}.
\begin{cooltheorem}\label{GEHA}
  Suppose $(M, M')$ is a pair of $C_p$-unital magma structures on the same coefficient system satisfying suitable interchange relations.
  Then, $M \simeq M'$ and each underlie a semi-Mackey functor;
  in particular, if the multiplications on $M^e$ and $M^{C_p}$ are invertible, then $M$ and $M'$  are isomorphic Mackey functors.
\end{cooltheorem}
For instance, given $X$ a $C_2$-unital magma valued in spaces, the $C_2$-space $UX$ defined by $UX^H = X(H)$ has an induced $C_2$-unital magma structure on $\upi_{\rho}(X)$ which interchanges with that of \cref{Homotopy coefficient system}.
\cref{GEHA} implies that these two structures agree and lift to a Mackey functor. 
\begin{example}
  The above argument confirms that the Mackey structure from Real Bott periodicity \cite{Atiyah} and the additive Mackey structure on Real vector bundles induce the same structure on $\upi_{2\rho} \BUR \simeq \upi_0 \BUR$.
\end{example}

To prove \cref{GEHA}, we embed $C_p$-unital magmas in the theory of \emph{algebras over $G$-operads} in the sense of \cite{Nardin} for $G$ a finite group;
in particular, we show in \cref{A2 section} that $C_p$-unital magmas are algebras over a particular $C_p$-operad $\AA_{2,C_p}^{\otimes}$ in $C_p$-coefficient systems valued in $\cC$, and spell out the correct interchange relations there.
We recommend that the reader familiarizes themself with the language of equivariant higher algebra via the introductions to \cite{EBV,Tensor}.

Crucially, the \emph{Boardman-Vogt tensor product}  of \cite{EBV} corepresents interchanging $G$-operad algebras:
\[
  \Alg_{\cO \otimes \cP}(\cD) \simeq \Alg_{\cO} \uAlg_{\cP}^{\otimes}(\cD).
\]
In particular, pairs of interchanging $C_p$-unital magma structures correspond with $\AA_{2,C_p}^{\otimes} \obv \AA_{2,C_p}^{\otimes}$-algebras.

Now, $G$-operads are manifestly homotopy-theoretic gadgets;
indeed, their algebras subsume the homotopy-coherent incomplete (semi-) Mackey functors of \cite{Blumberg_incomplete,Glasman,Cnossen_semiadditive} by \cite{Tensor,Marc}, the homotopy-coherent bi-incomplete Tambara functors of \cite{Blumberg-Bi-incomplete,Elmanto} by \cite{Chan,Cnossen_tambara}, and the algebraic structure on equivariant iterated loop spaces and their Thom spectra by \cite{Guillou,Horev_poincare}.
In particular, the first and second examples are incarnated by the \emph{weak $\cN_\infty$-operads} of \cite{EBV}, which are characterized by the fact that their nonempty structure spaces are contractible, and classified by the \emph{weak indexing category} given by the wide subcategory of maps of finite $G$-sets over which they supply structure:
\[
  A\cO \deq \cbr{T \rightarrow S \;\; \middle| \;\; \cO\text{-algebras are supplied maps } X^{\otimes T} \rightarrow X^{\otimes S}} \subset \FF_{G}.
\]
See \cite{Windex} for an overview of weak indexing categories.
In particular, given a subgroup $H \subset G$ and a finite $H$-set $S$, a $G$-operad $\cO$ has an \emph{$S$-ary structure space} $\cO(S)$ with action
\[
  \cO(S) \longrightarrow \Map\prn{X^{\otimes \Ind_H^G S} ,\; X^{\otimes G/H}},
\]
so $A\cO$ can be rigorously defined via vanishing of the structure spaces of $\cO$.

Conveniently, $\cO$-algebras in a $G$-symmetric monoidal $n$-category are canonically equivalent to algebras over the \emph{homotopy $n$-operad} $h_n \cO^{\otimes}$, whose structure spaces are $(n-1)$-truncations of the structure spaces of $\cO^{\otimes}$ \cite{EBV}.\footnote{Throughout this article, \emph{$n$-category} will be used to refer to \emph{$(n,1)$-categories}, i.e. $\infty$-categories whose mapping spaces are all $(n-1)$-truncated. The reader should feel free to think mostly in terms of familiar examples, such as the $n$-categoy of $(n-1)$-truncated spaces, of $(n-1)$-truncated connective spectra, of small $(n-1)$-categories, or of the hammock localization of chain complexes with homology concentrated in degrees $[d,d+n-1]$ for some uniform $d$.}
In particular, the structure spaces of $\cO^{\otimes}$ are $n$-connected if and only if $h_n \cO^{\otimes}$ possesses a (unique) equivalence with a weak $\cN_\infty$-operad;
then, $\cO$-algebras in coefficient systems valued in an $n$-category $\cD$ are (homotopy-coherent) incomplete semi-Mackey functors.
In this situation, we say that $\cO^{\otimes}$ is \emph{$n$-connected}.

From this, we identify \cref{GEHA} with the statement that $\AA_{2,C_p}^{\otimes} \obv \AA_{2,C_p}^{\otimes}$ is (0-)connected, together with the easy observation that $A\AA_{2,C_p} = \FF_{C_p}$,
so the corresponding incomplete semi-Mackey functors have all transfers.
Thus it suffices to prove the following equivariant generalization of \SY{Thm.}{1.0.1}{3}.
\begin{cooltheorem}[Equivariant Eckmann-Hilton argument]\label{HEHA}
  If $\cO^{\otimes}$ and $\cP^{\otimes}$ are $k$ and $\ell$-connected almost-unital $G$-operads with $A\cO = A\cP$, then $\cO^{\otimes} \otimes \cP^{\otimes}$ is $(k + \ell + 2)$-connected.
\end{cooltheorem}
All nonempty $G$-operads are $(-1)$-connected, so this extends \cref{GEHA} to equivariant higher algebra.
\begin{corollary}[Equivariant stabilization hypotheses]\label{Connectivity corollary}
  If $\cO^{\otimes}$ is a nonempty almost-unital $G$-operad, then $\cO^{\otimes (n+1)}$ is $(n-1)$-connected;
  in particular, for any $G$-symmetric monoidal $n$-category $\cC^{\otimes}$,
  \[
      U\colon \uCAlg^{\otimes}_{A\cO}(\cC) \xrightarrow{\;\;\;\;\; \sim \;\;\;\;\;} \overbrace{\uAlg^{\otimes}_{\cO} \cdots \uAlg^{\otimes}_{\cO}}^{(n+1)\text{-fold}}(\cC),
  \]
  where the $(n+1)$-fold tensor product is taken as a colimit (and composition as limit) 
  in the case $n = \infty$.
\end{corollary}
Here, $\CAlg_{I}$ refers to algebras over the weak $\cN_\infty$-operad associated with $I$.
For instance, this immediately implies a version of the \emph{Baez-Dolan stabilization hypothesis}, in this case stating that $(n+2)$-tuply $\EE_V$-monoidal $n$-categories are equivalent to $\EE_{\infty V}$-monoidal $n$-categories, i.e. $AV$-semi-Mackey functors valued in $n$-categories.

For another example, \cref{Connectivity corollary}, the lax $G$-symmetric monoidality of $\upi_0\colon \uSp^{\otimes}_G \rightarrow \underline{\mathrm{Mack}}^{\Box}_G(\Ab)$, and the results of \cite{Chan} or \cite{Cnossen_tambara} together construct a natural $A\cO$-Tambara structure on the 0th homotopy groups of $\cO \obv \cO$-ring $G$-spectra;\footnote{To construct this lax symmetric monoidality, first note that $\uSp_{G, \geq 0}^{\otimes} \subset \uSp_G^{\otimes}$ is closed under tensor products, so the localization $G$-functor $\uSp_G \rightarrow \uSp_{G,\geq 0}^{\otimes}$ is given a lax $G$-symmetric monoidal structure by \cref{Doctrinal adjunction}. 
  Moreover, to construct a lax $G$-symmetric monoidal structure on $\tau_{\leq 0} = \pi_0\colon \uSp_{G,\geq 0} \rightarrow \uSp_G^{\heart} = \uMack_G(\Ab)$, in light of \cite{Nardin} we need only note that indexed tensor products take $\pi_0$-equivalences to $\pi_0$-equivalences and that the resulting structure agrees with the usual one on Mackey functors;
the former follows by the same fact for $G = e$, conservativity of $\prod_{(H) \subset G} \Phi^H$, and the geometric fixed point formulae of \cite{HHR}.}
this and a forthcoming equivariant Dunn additivity result will construct a natural $AV$-Tambara structure on the 0th homotopy Mackey functors of $\EE_{2V}$-ring $G$-spectra.

Now, we may remove the assumption $A\cO = A\cP$ in \cref{HEHA}, but we will need a more refined notion of connectivity.
In general, given a weak indexing category $I$, we say that $\cO^{\otimes}$ is \emph{$k$-connected at $I$} if, for all elements of the corresponding weak indexing system
\[
  T \in \FF_{I,H} \deq \cbr{S \in \FF_H \;\; \middle| \;\; \Ind_H^G S \rightarrow [G/H] \in I},
\]
the structure space $\cO(T)$ is $k$-connected.
We define the \emph{connectivity function}
\[
  \Conn_{\cO}\colon \wIndCat_{G}^{a\uni} \longrightarrow \ZZ  \cup \cbr{\infty}
\]
by the formula $\Conn_{\cO}(I) \deq \mathrm{inf}\cbr{k \mid \cO^{\otimes} \text{ is } k \text{-connected at } I}$.
This is a $G$-operadic version of the \emph{connectivity dimension function} of a $G$-space (c.f. \cite[Def~1.1.(vi)]{Lewis_hurewicz}).

Now, 
$\prn{\ZZ \cup \cbr{\infty}}^{\wIndCat_G^{a\uni}}$ forms a \emph{pointwise} commutative monoidal poset, i.e.
\[
  f \leq g \;\;\;\; \iff \;\;\;\; \forall I, \; f(I) \leq g(I); \hspace{40pt} f + g(I) = f(I) + g(I).
\]
In this language, we will prove the following strengthening of \cref{HEHA}.
\begin{cooltheorem}\label{LEHA}
  Given $\cO^{\otimes},\cP^{\otimes}$ a pair of almost-unital $G$-operads, we have $\Conn_{\cO} + \Conn_{\cP} + 2 \leq \Conn_{\cO \otimes \cP}.$
\end{cooltheorem}
To put \cref{HEHA,LEHA} into context, note that a $G$-operad $\cO^{\otimes}$ is a weak $\cN_\infty$-operad if and only if $\Conn_{\cO}$ has all values $-2$ or $\infty$;
in this case, \cref{LEHA} says that weak $\cN_\infty$-operads are closed under tensor products and $\cN_{I \infty}^{\otimes} \obv \cN_{J \infty}^{\otimes}$ is classified by a weak indexing category contained in the join $I \vee J$.

Fortunately, this case of \cref{LEHA} is the difficult part of the main theorem of \cite{Tensor}.
To explain the main strategy of this, we must introduce a few definitions:
for $S \in \FF_H$, the \emph{$S$-indexed Wirthm\"uller map} in a (suitably pointed) $G$-$\infty$-category is defined to be the $S$-indexed semiadditive norm map as in \cite{Nardin-Stable,Cnossen_semiadditive};
that is, the $[H/K]$-indexed Wirthm\"uller map $W_{[H/K],X}\colon \Ind_K^H X \rightarrow \CoInd_K^H X$ is adjunct to the map 
\[
  X \longrightarrow \Res_K^H \CoInd_K^H X \simeq \prod_{g \in [K\backslash H / K]} \CoInd_{H \cap gKg^{-1}}^H \Res^H_{H \cap gKg^{-1}} X
\]
whose projection onto the factor indexed by the identity double coset is the identity and whose other projections are zero, and
the $\coprod_i [H/K_i]$-indexed Wirthm\"uller map 
\[
  W_{\coprod_i [H/K_i], (X_i)}\colon \coprod_{K_i}^H X_i  \simeq \coprod_i \Ind_{K_i}^H X_i \longrightarrow \prod_i \CoInd_{K_i}^H X_i \simeq \prod_{K_i}^H X_i
\]
is classified by the diagonal matrix whose $i$th entry is $W_{[H/K_i],X_i}$.

Key to \cite{Tensor} was the result that $\cO$-monoid spaces have $I$-indexed Wirthm\"uller isomorphisms 
if and only if $\cO^{\otimes}$ is $\infty$-connected at $I$.
In particular, this identified $\cN_{I\infty}^{\otimes}$ as the \emph{unique} $G$-operad $\cO^{\otimes}$ with $A\cO \leq I$ such that $\uMon_{\cO}(\cS)$ has $I$-indexed Wirthm\"uller isomorphisms.

Analogously, 
the key to this article will be to find \emph{$\ell$-connectivity of $\cO^\otimes$ at $I$} within the homotopy theory of $\cO$-monoids.
Explicitly, 
we say that a morphism $g\colon X \rightarrow Y$ in an $\infty$-category $\cC$ is \emph{$\ell$-truncated} if, for all $Z \in \cC$, the map of spaces $\Map(Z,X) \rightarrow \Map(Z,Y)$ is $\ell$-truncated, and $f\colon A \rightarrow B$ is \emph{$\ell$-connected} if, for all diagrams
  \[\begin{tikzcd}[ampersand replacement=\&]
	A \& X \\
	B \& Y
	\arrow[from=1-1, to=1-2]
	\arrow["f"', from=1-1, to=2-1]
	\arrow["g", from=1-2, to=2-2]
	\arrow[dashed, from=2-1, to=1-2, "h" description]
	\arrow[from=2-1, to=2-2]
\end{tikzcd}\]
such that $g$ is $\ell$-truncated, the space of lifts $h$ is contractible.
We will prove the following.
\begin{cooltheorem}\label{k-connected theorem}
  Let $\cP^{\otimes}$ be a $G$-operad, $I$ an almost unital weak indexing category, and $\ell$ a natural number.
  Then, the following conditions are equivalent:
  \begin{enumerate}[label={(\alph*)}]
    \item \label[condition]{cond:ell-connected at I} $\cP^{\otimes}$ is $\ell$-connected at $I$.
    \item \label[condition]{cond:topos wirthmuller maps} For all $n$-toposes $\cC$ (with $n \leq \infty$), 
      $I$-admissible $H$-sets $S \in \FF_{I,H}$, and
      $S$-tuples of $\cP$-monoids $(X_K) \in \prod_{[H/K] \in \Orb(S)}\Mon_{\Res_K^{G} \cP}\prn{\cC}$, 
      the $S$-indexed $\cP$-monoid Wirthm\"uller map
      \[
        W_{S,(X_K)}\colon \coprod_K^S X_K \longrightarrow \prod_K^S X_K
      \]
      is $\ell$-connected.
    \item \label[condition]{cond:spaces wirthmuller maps} For all 
      $I$-admissible $H$-sets $S \in \FF_{I,H}$ and
      $S$-tuples of $\cP$-monoids
      $(X_K) \in \prod_{[H/K] \in \Orb(S)}\Mon_{\Res_K^{G} \cP}(\cS)$, 
      the $S$-indexed $\cP$-monoid space Wirthm\"uller map
      \[
        W_{S, (X_K)}\colon \coprod_K^S X_K \longrightarrow \prod_K^S X_K
      \]
      is $\ell$-connected.
\end{enumerate}
\end{cooltheorem}
\begin{remark}\label{Topos-theoretic connectedness}
  \begin{enumerate}
    \item When $G = e$, only the implication (a) $\implies$ (b) is argued in \cite{Schlank}, but (b) $\implies$ (c) is obvious and, other than a diagram chase, the core idea of our argument that (c) $\implies$ (a) is already present in \cite{Schlank}.
    \item In the case that $\cC$ is an $n$-topos for some $0 \leq n \leq \infty$, the above definitions are equivalent to \emph{$\ell$-truncatedness} and \emph{$(\ell-1)$-connectiveness} in the sense of \httdef{6.5.1.10} by \SY{Lem.}{4.2.6}{29}
    and \httpropprop{6.5.1.12}{6.5.1.19}.
    \item While proving \cref{k-connected theorem}, we will verify that \cref{cond:topos wirthmuller maps} is further equivalent to the condition that the $\CoFr^H \cC$-map underlying $W_{S, (X_K)}$ is pointwise $\ell$-connected;
      moreover, \cref{cond:spaces wirthmuller maps} is equivalent to the condition that the underlying $H$-space map is $\ell$-connected, i.e. its associated maps on $J$-fixed point spaces are surjective on path components with $\ell$-connected fiber for each $J \subset H$.\qedhere
  \end{enumerate}
\end{remark}

\subsection*{Corollaries}
We highlight two main corollaries to the above results;
first, by combining them with the work done in \cite{Tensor}, we are able to characterize all $\obv$-idempotent almost-unital operads, or equivalently, the smashing localizations on $\Op_G^{a\uni}$.
To do so, define the full subcategory\def\Wirth{\mathrm{Wirth}}
\[
  \Op_{G}^{I-\Wirth} \deq \cbr{\cO^{\otimes} \;\; \middle| \;\; \forall \; S \in \FF_{I,H} \text{ and }  \cC^{\otimes} \in \Cat_{G}^{\otimes}, \;\; \text{ we have} \;\;\bigotimes^S \simeq \coprod^S \; \text{ in } \uAlg_{\cO}(\cC)} \subset \Op_{G}^{a\uni}.
\]
We will conclude that $\cN_\infty$-operads and indexing systems are additionally equivalent to \emph{smashing localizations} on a suitable category of $G$-operads:
\begin{coolcorollary}\label{cor:smashing}
  The following commutative diagram consists of isomorphisms of posets
  \[
    \begin{tikzcd}[ampersand replacement=\&]
      {    \wIndSys^{a\uni}_{G} } \& {\cbr{\text{\rm Subterminal objects of $\Op_G^{a\uni}$}}} \\
	\& {\cbr{\text{\rm Smashing localizations of } \Op_{G}^{a\uni} \text{ \rm under reverse inclusion}}}
  \arrow["{\cN_{\bullet\infty}^{\otimes}}"', shift right, from=1-1, to=1-2]
  \arrow["{\Op_{G}^{\bullet-\Wirth}}"', end anchor={west}, from=1-1, to=2-2]
	\arrow["A"', shift right, from=1-2, to=1-1]
	\arrow["{\Mod(-)}", shift left=2, from=1-2, to=2-2]
	\arrow["{L \EE_{0}^{\otimes}}", from=2-2, to=1-2]
\end{tikzcd}
  \]
\end{coolcorollary}
A striking corollary of this is that there are finitely many smashing localizations on $\Op_{G}^{a\uni}$;
moreover, they have rich combinatorial structure, as they are naturally cocartesian-fibered over transfer systems, giving e.g. a cocartesian fibration from smashing localizations on $\Op_{C_{p^n}}^{a\uni}$ to the $(n+2)$nd associahedron, whose fibers can be explicitly described \cite{Balchin,Windex}.

Before moving on, we point out another corollary of \cref{LEHA}.
Given $I$ an indexing category, let $\Sp_I$ be the $\infty$-category presented by Blumberg-Hill's stable model category of $I$-spectra \cite{Blumberg_stable}.
Let $\uAlg_{\cO}(\cS_{G, [k,\ell]}) \subset \uAlg_{\cO}(\cS_G)$ be the full subcategory spanned by $\cO$-monoid $G$-spaces whose underlying $G$-spaces have $\ZZ$-graded coefficient system homotopy groups concentrated in degrees $[k,\ell]$, and similarly define $\Sp_{G, [k,\ell]}$.
We will prove the following.
\begin{coolcorollary}\label{I-spectrum coolcorollary}
  Fix $\cO^{\otimes}$ a reduced $G$-operad with $\cO(2 \cdot *_G) \neq 0$, and  $0 \leq k \leq \ell \leq \infty$ numbers.
  There is an equivalence
  \[
    \begin{tikzcd}[column sep=small, ampersand replacement=\&]
	{\Sp_{G, [k,\ell]}} \& {\overbrace{\Alg_{\cO} \uAlg_{\cO}^{\otimes} \cdots \uAlg_{\cO}^{\otimes}}^{(\ell-k+2)\text{-fold}}(\ucS_{G, [k,\ell]})} \\
	{\cS_{G, [k,\ell]}}
	\arrow["\simeq"{description}, draw=none, from=1-1, to=1-2]
	\arrow["{\Omega^{\infty}}"', from=1-1, to=2-1]
	\arrow["U", from=1-2, to=2-1]
\end{tikzcd}
  \]
  That is, to equip a $(k-1)$-connected and $\ell$-truncated $G$-space $X$ with $(\ell - k+2)$-many interchanging $\cO$-algebra structures is equivalent to realizing $X$ as the 0th $G$-space of an $A\cO$-spectrum.
\end{coolcorollary}
\begin{remark}
  Qualitatively, this is much weaker than the result we find nonequivariantly, e.g. from \SY{Thm.}{5.2.2}{37}.
  This is because norms \emph{need not exacerbate connectivity};
  indeed, given $X \in \cS_{H,*}$ such that $X^H$ is \emph{not} $n$-connected, $\prn{\CoInd_H^G X}^G \simeq X^H$ is not $n$-connected, so $\CoInd_H^G X$ is not $n$-connected.
  
  The corresponding fact in topology is that, if $A$ is a coefficient system of Abelian groups which does not extend to an $I$-Mackey functor, then the Eilenberg-Mac Lane $G$-space $K(A,n)$ is not the $n$th space of any $I$-spectrum.
  The author suspects that a version of the strong result follows for concentration in particular \emph{regular slice degrees}, but we will not discuss that here.
\end{remark}

To construct an infinite loop space theory for $I$-spectra, one is left with the following question.
\begin{question}\label{Loop space theory question}
  Given an indexing category $I$, does there exist a reduced $G$-operad $\cO^{\otimes}$ with $A\cO = I$ and a space $S^I$ such that $\cO$-monoid structures on a connected $G$-space $X$ are equivalent to $S^I$-loop space structures?
\end{question}

\begin{remark}
  \cref{I-spectrum coolcorollary} has a strong philosophical implication running transverse to \cref{Loop space theory question}:
  regardless of the topology, it constructs a flexible machine which inputs unital equivariant algebraic theories and outputs towers of $\infty$-categories converging to equivariant stable homotopy theory.
  For instance, iterating algebras over Rubin's free or associative $N$-opeards \cite{Rubin} yields such a tower converging to arbitrary $\Sp_I$.

  In essence, \cref{I-spectrum coolcorollary} \emph{takes the loops out of equivariant infinite loop space theory}, extending algebraic versions of the theory to arbitrary incompletely stable categories regardless of the answer to \cref{Loop space theory question}.
\end{remark}  

\subsection*{Philosophy}
The following significantly motivated this article and its prequels \cite{Windex,EBV,Tensor}.
\begin{question}
  For what higher-algebraic, universal, and operadic reasons do $\cN_\infty$ operads arise?
\end{question}
Of course, there are preexisting higher-algebraic reasons:
the several incomplete variants of the \emph{spectral Mackey functor theorem} verify that $\cN_{I\infty}$-monoids are intimately connected with $I$-admissible Wirthm\"uller isomorphisms, which are close to the heart of equivariant stable homotopy theory \cite{Guillou,Nardin-Stable,Cnossen_semiadditive,Blumberg_incomplete,Marc}.
Since the value $\uAlg_{\cN_{I\infty}}(\ucS_{G})$ uniquely pins $\cN_{I\infty}^{\otimes}$ as a $G$-operad \cite{Tensor}, this characterizes $\cN_\infty$-operads.
This is not a complete answer, as it requires us to care about $I$-indexed Wirthm\"uller isomorphisms a-priori;
that is, while our reason is higher-algebraic and universal, it is not quite operadic in its philosophy.

Now, if we admit \emph{weak} $\cN_\infty$-operads, a universal operadic characterization is easy to come by: they are precisely the subterminal objects of $\Op_G$.
Unfortunately, \emph{algebra} lives in the mapping spaces from one-object $G$-operads to $G$-symmetric monoidal $\infty$-categories, and no nontrivial $\cN_\infty$-operads are $G$-symmetric monoidal $\infty$-categories, so the \emph{mapping-in} property identifying weak $\cN_\infty$-operads is not higher-algebraic in nature.

In the almost-unital locus (or, for that matter, the unital locus), \cref{Classification of idempotents} gives a characterization with all three properties:
almost-unital weak $\cN_\infty$-operads are characterized universally as the corepresenting $G$-operads at the limit of of (infinitary) Eckmann-Hilton arguments in equivariant higher algebra.
In essence, they are exactly the class of $G$-operads $\cO^{\otimes}$ possessing the property that the natural forgetful functor
\[
  \Alg_{\cO} \uAlg_{\cO}^{\otimes}(\ucS_G) \rightarrow \Alg_{\cO}(\ucS_G)
\]
forgetting the inner one $\cO$-algebra structure is an equivalence.

\subsection*{The strategy}
In \cref{Reduction subsection} we reduce \cref{LEHA,HEHA} to the case of \cref{LEHA} where $\cO,\cP$ are \emph{unital}.
In this case, we perform a similar reduction to \cite{Schlank};
namely, by examining the free $\cO$-algebra monad, we reduce this to \emph{$(k + 1)$-connectivity of the reduced endomorphism $A\cO$-operad in $\uMon_{\cP}(\cC)^{I-\times}$} in the case $\cC$ is the $G$-$\infty$-category of coefficient systems in a presheaf $\infty$-topos.

We express the structure space $\End^{\red}_X\prn{\uMon_{\cO}(\cC)^{I-\times}}(S)$ as the spaces of lifts of $\Delta\colon X^{\sqcup S} \rightarrow X$ along the $S$-indexed Wirthm\"uller map $W_{X,S}\colon X^{\sqcup S} \rightarrow X^{\times S}$, which is directly related to truncatedness of $X$ and connectedness of $W_{X,S}$;
hence it suffices to prove \cref{k-connected theorem} in the unital case.

We finish by directly relating $\ell$-connectivity of $W_{X,S}$ in $\Mon_{\cO}(\cC)$ and $\Mon_{\cO}(\tau_{\leq \ell} \cC)$, reducing \cref{k-connected theorem} to the fact that $\Mon_{\cO}(\tau_{\leq \ell} \cC)$ is $I$-semiadditive when $\cO$ is $\ell$-connected at $I$, which we verified in \cite{Tensor}.

\subsection*{Acknowledgements}
This article is greatly influenced by the work of Schlank-Yanovski \cite{Schlank}, which recovers almost all of the results of this article in the case that $G$ is the trivial group, and has additionally been influential to my thinking in the previous articles \cite{EBV,Tensor}.
In general, I'd like to thank my advisor Mike Hopkins for several helpful conversations on this material and Branko Juran for comments on an early draft of this article.

\resumetocwriting
\section{Preliminaries} 
The rest of this paper replaces the orbit category $\cO_G$ with an arbitrary atomic orbital $\infty$-category $\cT$;
that is, we assume that all retracts in $\cT$ are equivalences and that the finite coproduct completion $\FF_{\cT} \deq \cT^{\sqcup}$ has pullbacks.
We will prove \cref{HEHA,LEHA,k-connected theorem} in that level of generality.
We encourage the reader to either globally specialize to $\cT = \cO_G$ or familiarize themself with the atomic orbital setting via \cite{EBV} .

We begin in \cref{sub:windex prelims,sub:CMon prelims} by recalling the simultaneous generalization and weakening of Blumberg-Hill's $G$-indexing systems and $I$-Mackey functors to $\cT$-\emph{weak} indexing systems and $I$-commutative monoids.
We go on to \cref{sub:op prelims} where we recall the relevant background from \cite{Nardin,EBV,Tensor} on $\cT$-operads.

\subsection{Preliminaries on \tcT-\tinfty-categories and weak indexing systems}\label{sub:windex prelims}
Recall that a \emph{$\cT$-coefficient system} is a functor out of $\cT^{\op}$:
\[
  \CoFr^{\cT}(\cC) \deq \Fun(\cT^{\op},\cC).
\]
Extending Elmendorf's theorem, we define $d$-truncated $\cT$-spaces and $\cT$-$d$-categories as coefficient systems: 
\[
  \cS_{\cT, \leq d} \deq \CoFr^{\cT}(\cS_{\leq d}); \hspace{50pt} \Cat_{\cT, d} \deq \CoFr^{\cT}(\Cat_d).
\]
We write $\Cat_{\cT} \deq \Cat_{\cT,\infty}$ and $\cS_{\cT} \deq \cS_{\cT, \leq \infty}$.
Given a $\cT$-$\infty$-category $\cC$, we write $\cC_V$ for the value $\cC(V)$ and $\Res_V^W\colon \cC_W \rightarrow \cC_V$ for the functoriality under a map  $V \rightarrow W$.
The $\infty$-category of $\cT$-coefficient systems lifts to a $\cT$-$\infty$-category with $V$-value the $\cT_{/V}$-coefficient systems
\[
  \uCoFr^{\cT}(\cC)_V \deq \CoFr^{\cT_{/V}}(\cC);
\]
the functoriality is given by restriction.
We acquire $\cT$-$\infty$-categories $\ucS_{\cT, \leq d}$ and $\uCat_{\cT,d}$ similarly.
\begin{example}
  We may define a $\cT$-$\infty$-category by $\uFF_{\cT}$ by values
  \[
    \prn{\uFF_{\cT}}_V \deq \FF_{\cT,/V} \simeq \FF_{\cT_{/V}}
  \]
  with functoriality given by pullback.
  We write $\FF_V \deq \FF_{\cT, /V}$.
  This is a $\cT$-$1$-category since $\cT_{/V}$ is a 1-category \NSprop{2.5.1}.
\end{example}
\begin{example}
  Given $\cC$ an arbitrary $n$-category, $\uCoFr^{\cT}(\cC)$ is a $\cT$-$n$-category \httcor{2.3.4.8}.
  In particular, if $\cC$ is an $\infty$-topos and $\tau_{\leq n-1} \cC$ its $n$-topos of $(n-1)$-truncated objects, then $\uCoFr^{\cT}(\tau_{\leq n-1} \cC)$ is a $\cT$-$n$-category.
\end{example}
\begin{example}[{\cite[Thm.~\href{https://arxiv.org/pdf/1608.03657v1\#nul.9.7}{9.7}]{Barwick-Parameterized}}]
  The $\infty$-category of $\cT$-$\infty$-categories is Cartesian closed with internal hom characterized by values
  \[
    \uFun_{\cT}(\cC,\cD)_V \simeq \Fun_{\cT_{/V}}(\Res_V^{\cT} \cC, \Res_V^{\cT} \cD),
  \]
  where $\Res_V^{\cT}\colon \Cat_{\cT} \rightarrow \Cat_{\cT_{/V}}$ is pullback and $\Fun_{\cT_{/V}}(-,-)$ denotes the $\infty$-category of functors of cocartesian fibrations between unstraigthenings over $\prn{\cT_{/V}}^{\op}$. 
\end{example}
\begin{example}\label{Inflation example}
  We refer to the adjunction between limits and constant diagrams as the \emph{inflation and fixed point} adjunction
    \[
      \begin{tikzcd}[ampersand replacement=\&]
        \Cat \& {\Cat_{\cT}}
        \arrow[""{name=0, anchor=center, inner sep=0}, "{\Infl_e^{\cT}}", curve={height=-12pt}, from=1-1, to=1-2]
        \arrow[""{name=1, anchor=center, inner sep=0}, "{\Gamma^{\cT}}", curve={height=-12pt}, from=1-2, to=1-1]
        \arrow["\dashv"{anchor=center, rotate=-90}, draw=none, from=0, to=1]
      \end{tikzcd}
    \]
  In the case that $\cT$ has a terminal object $V$, the image of $\Infl_e^{\cT}$ consists of the $\cT$-$\infty$-categories whose restriction functors $\Res^V_W$ are all equivalences.
  In any case, we may string together natural equivalences
  \begin{align*}
    \uFun_{\cT}\prn{\Infl_e^{\cT} K, \uCoFr^{\cT} \cC}_V
    &\simeq \Fun_V\prn{\Infl_e^{\cT_{/V}} K, \uCoFr^{\cT_{/V}}\cC}\\
    &\simeq \Fun\prn{K, \Fun\prn{\prn{\cT_{/V}}^{\op}, \cC}}\\
    &\simeq \Fun\prn{\prn{\cT_{/V}}^{\op}, \Fun(K, \cC)}\\
    &\simeq \uCoFr^{\cT}\prn{\cC^K}_V
  \end{align*}
  to construct a $\cT$-equivalence $\uFun_{\cT}\prn{\Infl_e^{\cT} K, \uCoFr^{\cT} \cC} \simeq \uCoFr^{\cT}\prn{\cC^K}$;
  in particular, choosing $\cC = \cK$, $\cT$-coefficient systems in presheaves of spaces on $K$ can equivalently be realized as $\cT$-equivariant presheaves of $\cT$-spaces on $K$ with trivial $\cT$-equivariant structure.
  We henceforth write
  \[
    \ucS_{\cT, \leq n}^K \deq \uCoFr^{\cT}\prn{\cS_{\leq n}^K}; \hspace{50pt}
    \ucS_{\cT}^K \deq \uCoFr^{\cT}\prn{\cS^K}.\qedhere
  \]
\end{example}

Given $V \in \cT$ an orbit and $S \in \FF_V$ a finite $V$-set, we write $\varphi_{SV}\colon \Ind_V^{\cT} S \rightarrow V$ for the corresponding map in $\FF_{\cT}$, and we write
\[
  \cC_S \deq \prod_{U \in \Orb(S)} \cC_U \simeq \Fun_{\cT}\prn{\Ind_V^{\cT} S, \cC}.
\]
Pullback along the structure map $\varphi_{SV}$ yields an \emph{indexed diagonal} functor
\[
  \Delta^S\colon \cC_V \rightarrow \cC_S;
\]
its values are $\Delta^S X = (\Res_U^V X)_{U \in \Orb(S)}$.
The \emph{$S$-indexed coproduct} (if it exists) is the left adjoint $\coprod^S\colon \cC_S \rightarrow \cC_V$ to $\Delta^S$, and the \emph{$S$-indexed product} $\prod^S\colon \cC_S \rightarrow \cC_V$ is the right adjoint.
\begin{notation}\label{Distinguished fixed point}
  In the case $U \rightarrow V$ is a map of orbits, considered as an element of $\FF_V = \FF_{\cT, /V}$, we write
  \[
    \Ind_U^V(-) \deq \coprod^{U \rightarrow V}(-); \hspace{30pt} \CoInd_U^V(-) \deq \prod^{U \rightarrow V}(-),
  \]
  so that $\Ind_U^V \dashv \Res_U^V \dashv \CoInd_U^V$. We refer to these as \emph{induction and coinduction}.
  Note that
    \[
      \Ind_U^V\colon \FF_{\cT, /U} \rightarrow \FF_{\cT, /V}
    \]
  is postcomposition.
  We call $S^V = \Hom_{\FF_V}(*_V, S)$ the \emph{fixed points} and define the \emph{distinguished fixed point} 
  \[\begin{tikzcd}[ampersand replacement=\&]
    U = \Ind_U^{\cT} *_U \\
	\& {\Ind_U^{\cT} \Res_U^V \Ind_U^{V} *_V} \& {U = \Ind_V^\cT \Ind_U^V*_U} \\
	\& U \& V
  \arrow["{\Ind_U^{\cT} \delta}"{description}, dashed, from=1-1, to=2-2]
	\arrow[curve={height=-12pt}, equals, from=1-1, to=2-3]
	\arrow[curve={height=18pt}, equals, from=1-1, to=3-2]
	\arrow[from=2-2, to=2-3]
	\arrow[from=2-2, to=3-2]
	\arrow["\lrcorner"{anchor=center, pos=0.125}, draw=none, from=2-2, to=3-3]
	\arrow[from=2-3, to=3-3]
	\arrow[from=3-2, to=3-3]
\end{tikzcd}\]
 Note that, since $\cT$ is atomic, $\delta\colon *_U \rightarrow \Res_U^V \Ind_U^V *_U$ is a summand inclusion.
 In analogy to equivariant homotopy theory, we suggest the reader view $\delta_U$ as ``the identity coset fixed point.''
\end{notation}
These are the ur-examples of \emph{equivariantly indexed operations}, whose combinatorics we control using \emph{weak indexing systems}.
\begin{definition}
  A \emph{one-color weak indexing system} is a full $\cT$-subcategory $\uFF_{I} \subset \uFF_{\cT}$ which is closed under $\uFF_I$-indexed coproducts and contains $*_{V}$ for all $V \in \cT$.
  A \emph{one-color weak indexing category} is a pullback-stable wide subcategory $I \subset \FF_{\cT}$ subject to the condition that $\coprod_i\prn{T_i \rightarrow S_i}$ lies in $I$ if and only if each map $T_i \rightarrow S_i$ lies in $I$.
\end{definition}

Given $I$ a one-color weak indexing category, we define the \emph{$I$-admissible $V$-sets} as
\[
  \uFF_I \deq \cbr{S \mid \Ind_V^{\cT} S \rightarrow V \in I} \subset \uFF_{\cT};
\]
we verified in \cite[Thm.~\href{https://arxiv.org/pdf/2409.01377v2\#cooltheorem.1}{A}]{Windex} that $\uFF_{(-)}$ furnishes an equivalence between one-color weak indexing systems and one-color weak indexing categories, so we conflate these notions.
For the following examples, a full subcategory $\cF \subset \cT$ is called a \emph{$\cT$-family} if, whenever there exists a morphism $V \rightarrow W$ with $W \in \cF$, we have $V \in \cF$.
\begin{example}
  The terminal one-color weak indexing system is $\uFF_{\cT}$.
  Given a family $\cF \subset \cT$, we define
  \begin{align*}
    \prn{\uFF_{\triv}}_V 
        &\deq \cbr{*_V}\\
    \prn{\uFF_{0,\cF}}_V 
        &\deq \begin{cases}
              \cbr{\emptyset_V,*_V} & V \in \cF \\   
              \cbr{*_V} & \mathrm{otherwise}.
            \end{cases}\\
    \prn{\uFF_{\infty}}_V
        &\deq \cbr{n \cdot *_V \mid n \in \NN}.
  \end{align*}
  The corresponding one-color weak indexing categories are denoted $I_{\triv}, I_{0,\cF}, I_\infty$.
\end{example}

\begin{construction}[\Windlem{1.24}]
  Given $I$ a one-color weak indexing category, We define the \emph{unit family}
  \[
    \upsilon(I) \deq \cbr{V \in \cT \mid \emptyset_V \in \prn{\uFF_{I}}_V} \subset \cT.\qedhere
  \]
\end{construction}
We say that $\uFF_I$ is \emph{almost-unital} if, whenever $\cbr{*_V} \subsetneq \FF_{I,V}$, we have $\emptyset_V \in \FF_{I,V}$;
that is, $\uFF_I$ is unital over all orbits for which $\uFF_I$ has nontrivial arities.
We say $\uFF_I$ is \emph{unital} if $\emptyset_V \in \FF_{I,V}$ for all $V$, i.e. $\upsilon(I) = \cT$.
We denote the corresponding posets (under inclusion) as
\[
  \wIndSys_{\cT}^{\uni} \subset \wIndSys_{\cT}^{a\uni} \subset \wIndSys_{\cT}^{\oc}.
\]
The following lemma is an important aspect of almost-unitality.
\begin{lemma}[\Windlem{1.25}]\label{Almost-unital summands}
  $\uFF_I$ is almost-unital if and only if, for each \emph{nonempty} summand $S \subset T$ with $T \in \FF_{I,V}$ we have $S \in \FF_{I,V}$.
\end{lemma}

\subsection{Preliminaries on \tI-commutative monoids and \tI-symmetric monoidal \tinfty-categories}\label{sub:CMon prelims}
Let $I$ be a one-color weak indexing category.
The pair $(\FF_{\cT},I)$ is a \emph{span pair} in the sense of \eh{Def.}{2.2.1} (i.e. $(\FF_{\cT},I,I)$ is an \emph{adequate triple} in the sense of \cite[Def.~\href{https://arxiv.org/pdf/1404.0108v2\#page=15}{5.2}]{Barwick1}), so it yields a wide subcategory
\[
  \Span_I(\FF_{\cT}) \hookrightarrow \Span(\FF_{\cT})
\]
of the effective Burnside $\infty$-category with morphisms the spans $X \leftarrow R \xrightarrow{f} Y$ with $f \in I$.
Given $I$ a one-color weak indexing category and $\cC$ an $\infty$-category, we define the $\infty$-category of \emph{$I$-commutative monoids in $\cC$}
\[
  \CMon_I(\cC) \deq \Fun^{\times}(\Span_I(\FF_{\cT}), \cC).
\]
We define the $\infty$-category of \emph{small $I$-symmetric monoidal $\infty$-categories} as
\[
  \Cat_I^{\otimes} \deq \CMon_I(\Cat).
\]
We henceforth ignore size issues and omit the adjective ``small.''
Given an $I$-symmetric monoidal $\infty$-category $\cC$ and $S \in \FF_{I,V}$ an $I$-admissible $V$-set, we denote the functoriality of $\cC^{\otimes}$ under the structure map $\Ind_S^{\cT} S = \Ind_S^{\cT} S \rightarrow V$ by
\[
  \bigotimes^S\colon \cC_S \rightarrow \cC_V.
\]
Now, in the case that $I$ is almost-unital, the orbit collapse factorization $\Ind_V^{\cT} S \rightarrow \coprod_{U \in \Orb(S)} V \rightarrow V$ lies in $I$ by \cref{Almost-unital summands}, so we acquire a natural equivalence
\[
  \bigotimes^S_U X_U \simeq \bigotimes_{U \in \Orb(S)} N_U^V X_U,
\]
where we write $N_U^V \coloneqq \bigotimes^{U \rightarrow V}$.
A similar factorization holds for indexed (co)products.

If $I$ is almost-unital, $S \in \FF_{I,V}$ is $I$-admissible, and $1_U \in \cC_U$ is initial whenever it exists, then given an $S$-indexed tuple $(X_U) \in \cC_S$ in an $I$-symmetric monoidal $\infty$-category with $S$-indexed coproducts, we define an \emph{$S$-indexed tensor Wirthm\"uller map}
\[
  W_{S,(X_U)}\colon \coprod_{U}^S X_U \longrightarrow \bigotimes_U^S X_U
\]
uniquely extending the by composite maps $\Ind_W^V X_W \hookrightarrow \coprod_U^S X_U \rightarrow \bigotimes_U^S X_U$, which are adjunct to the map
\[
  \iota_W\colon X_W \simeq X_W \otimes \bigotimes^{\Res_W^V S - \delta(W)}_U 1_U \xrightarrow{\; \; \prn{\id, \eta} \; \;} X_W \otimes \bigotimes^{\Res_W^V S - \delta(W)}_U X_U \simeq \Res_W^V \bigotimes^{S}_U X_U;
\]
here, $\partial(W) \subset \Res_W^V \Ind_W^V *_W \subset \Res_W^V S$ is the distinguished fixed point of \cref{Distinguished fixed point}.
Intuitively, on the $W$'th factor, $W_{S,(X_U)}$ takes $x$ to the simple tensor with the Wirthm\"uller image of $x$ in the $W$'th place and units elsewhere.
Given $J \subset I$, we say that $\cC$ is \emph{$J$-cocartesian} if $W_{S,(X_U)}$ is an equivalence for all $S \in \uFF_K$ and $(X_U) \in \cC_S$, and we say that $\cC$ is \emph{$J$-cartesian} if its ``vertical opposite'' 
\[
  \Span_I(\FF_{\cT}) \xrightarrow{\cC^{\otimes}} \Cat \xrightarrow{\op} \Cat
\]
is a $J$-cocartesian $I$-symmetric monoidal $\infty$-category.

This is functorial.
To say how, let $\Cat_I^{\sqcup} \subset \Cat_{\cT}$ be the replete subcategory of $\cT$-categories with $I$-indexed coproducts and $I$-coproduct preserving $\cT$-functors, and similarly define $\Cat_I^{\times} \subset \Cat_{\cT}$.
\begin{proposition}[{\cite[Thm.~\href{https://arxiv.org/pdf/2504.02143v2\#cooltheorem.1}{A}]{Tensor} and \stecor{1.65}}]\label{cor:cart}
  The $I$-cocartesian and $I$-cartesian $I$-symmetric monoidal categories form the essential image of fully faithful inclusions
  \[
    \begin{tikzcd}[ampersand replacement=\&, row sep=small, column sep=huge]
      {\Cat_I^{\sqcup}} \& {\Cat_I^{\otimes}} \& {\Cat_I^{\times}} \\
      \& {\Cat_I}
      \arrow["{(-)^{I-\sqcup}}", hook', from=1-1, to=1-2]
      \arrow[hook, from=1-1, to=2-2]
      \arrow["U", from=1-2, to=2-2]
      \arrow[hook', from=1-3, to=2-2]
      \arrow["{(-)^{I-\sqcup}}"', hook, from=1-3, to=1-2]
    \end{tikzcd}
  \]
  Moreover, $\cC$ is $I$-semiadditive in the sense of
  \sta{Def.}{4.5.1} if and only if there exists an equivalence $\cC^{I-\times} \simeq \cC^{I-\sqcup}$ lying over the identity endofunctor of $\cC$.
\end{proposition}

\subsection{Naive preliminaries on \tI-operads}\label{sub:op prelims}
In \cite{Nardin}, an $\infty$-category $\Op_{\cT}$ of \emph{$\cT$-operads} was introduced, and in \cite{EBV,Tensor} it was given a symmetric monoidal closed $\cT$-$\infty$-category structure $\uOp_{\cT}^{\otimes}$.
We review the relevant formal properties here;
in particular, 
we will only use properties of $\uOp_{\cT}^{\otimes}$ and the various functors
\[\begin{tikzcd}[ampersand replacement=\&, column sep=huge]
	{\Cat_{\cT}^{\otimes}} \& {\Op_{\cT}} \& {\Fun(\Tot \uSigma_{\cT}, \cS)} \\
  {\Cat_{\cT}} \& \Cat_{\cT} \& \Cat_{\cT}
	\arrow[hook', from=1-1, to=1-2]
	\arrow["\sseq", from=1-2, to=1-3]
	\arrow["U"{description}, from=1-2, to=2-1]
	\arrow["{\uAlg_{(-)}(\cC)}"{description}, from=1-2, to=2-2]
	\arrow["{\uAlg_{\cP}(-)}"{description}, from=1-2, to=2-3]
\end{tikzcd}\]
In this way, this paper can be considered agnostic to the presentation of $\uOp_{\cT}^{\otimes}$ and the above functors.

\subsubsection{$\cT$-symmetric sequences and $I$-operads}
Writing $\uSigma_{\cT}$ for the composite $\cT$-$\infty$-category
\[
  \cT^{\op} \xrightarrow{\uFF_{\cT}} \Cat \xrightarrow{(-)^{\simeq}} \cS \hookrightarrow \Cat
\]
and writing $\Tot\colon \Cat_{\cT} \simeq \Cat_{/\cT^{\op}}^{\cocart} \rightarrow \Cat$ for the total category functor, in \EBVsubsec{2.3} we defined a \emph{underlying $\cT$-symmetric sequence} functor
\[
  \cO(-)\colon \Op_{\cT} \rightarrow \Fun(\Tot \uSigma_{\cT}, \cS).
\]
Moreover, in \ns{Def.}{2.1.7}, an \emph{underlying $\cT$-$\infty$-category} $U\colon \Op_{\cT} \rightarrow \Cat_{\cT}$ was defined.
We use this for the following definition.
\begin{definition}
  We say that an $I$-operad $\cO^{\otimes}$ \emph{has at least one color} if $\cO(*_V) \neq \emptyset$ for all $V \in \cT$ and \emph{has one color} if $\cO(*_V) \simeq *$ for all $V \in \cT$.
  We say that $\cO$ has \emph{connected colors} if $\pi_0 U\cO = *_{\cT}$, and we write $\Op_{\cT}^* \subset \Op_{\cT}$ for the full subcategory of $\cT$-operads with connected colors.
\end{definition}
\begin{proposition}[{\EBVsubsec{2.3}}]\label{thm:free operads}
  The functor $\cO(-)\colon \Op^*_{\cT} \rightarrow \Fun(\Tot \uSigma_{\cT}, \cS)$ has a left adjoint $\Fr$;
  in particular, letting $\Fr_{\Op}(S)$ be the free $\cT$-operad on the left Kan extended $\cT$-symmetric sequence
  \[\begin{tikzcd}[ampersand replacement=\&]
	{\cbr{S}} \& \cS \\
	{\Tot \uSigma_{\cT}},
	\arrow[hook, from=1-1, to=1-2]
	\arrow[hook, from=1-1, to=2-1]
	\arrow[""{name=0, anchor=center, inner sep=0}, "{\Fr_{\Sigma, S}(*)}"', dashed, from=2-1, to=1-2]
	\arrow[shorten >=3pt, Rightarrow, from=1-1, to=0]
\end{tikzcd}\]
  the adjunctions construct a natural equivalence
  \[
    \Alg_{\Fr_{\Op}(S)}(\cO) \simeq \cO(S).
  \] 
  Moreover, the restricted functor $\cO(-)\colon \Op_{\cT}^{*} \rightarrow \Fun(\Tot \uSigma_{\cT},\cS)$ is monadic.
\end{proposition}
In particular, identifying an object of $\Tot \uSigma_{\cT}$ with a pair $(V,S)$ where $V \in \cT$ and $S \in \FF_V$, $\cT$-operads are identified conservatively by the functor
\[
  \cO \mapsto \prod_{V,S} \cO(S).
\]
Intuitively, we view $\cO(S)$ as the space of \emph{$S$-ary operations} $\prn{\Res_V^{\cT} X}^{\otimes S} \rightarrow \Res_V^{\cT} X$ borne by an $\cO$-algebra $X$.
This technology allowed us to define the \emph{arity support} subcategory
\[
  A\cO \deq \cbr{T \rightarrow S \;\; \middle| \;\; \prod_{U \in \Orb(S)} \cO(T \times_S U) \neq \emptyset} \subset \FF_{\cT};
\]
which we verified in \EBVprop{2.88} to be a weak indexing category.
In fact, we verified in \EBVcor{2.91} that the essential surjection associated with $A$ possesses a fully faithful right adjoint
\[
  \begin{tikzcd}[ampersand replacement=\&]
    {\Op_{\cT}} \& {\wIndCat_{\cT}};
    \arrow[""{name=0, anchor=center, inner sep=0}, "A", curve={height=-12pt}, from=1-1, to=1-2]
    \arrow[""{name=1, anchor=center, inner sep=0}, "{\cN_{(-)\infty}^{\otimes}}", curve={height=-12pt}, hook', from=1-2, to=1-1]
    \arrow["\dashv"{anchor=center, rotate=-90}, draw=none, from=0, to=1]
  \end{tikzcd}
\]
we refer to the $\cT$-operad $\cN_{I\infty}^{\otimes}$ as the \emph{weak $\cN_\infty$-operad associated with $I$}.
Now, we further verified in \EBVcor{2.82} that, given a $\cT$-operad $\cO^{\otimes}$, the unique map $\cO^{\otimes} \rightarrow \Comm_{\cT}^{\otimes}$ is a monomorphism if and only if the counit map $\cO^{\otimes} \rightarrow \cN_{A\cO}^{\otimes}$ is an equivalence;
in particular, we acquire an equality of full subcategories
\[
  \Op_{\cT, /\cN^{\otimes}_{I \infty}} = A^{-1}(\wIndCat_{\cT, \leq I}) \;\; \subset \;\; \Op_{\cT},
\]
and a full subcategory of $\Op_{\cT}$ has a terminal object if and only if it is of this form.
We refer to $\Op_I \deq \Op_{\cT, /\cN_{I\infty}^{\otimes}}$ as the $\infty$-category of \emph{$I$-operads};
see \EBVprop{2.39} for an intrinsic characterization of $\Op_I$.

Monomorphisms are right-cancellable, so all inclusions $I \subset J$ induce monomorphisms $\iota_I^J\colon \cN_{I \infty}^{\otimes} \rightarrow \cN_{J \infty}^{\otimes}$;
in other words, the push-pull adjunction
\[
  \begin{tikzcd}[ampersand replacement=\&]
    {\Op_I} \& {\Op_J}
    \arrow[""{name=0, anchor=center, inner sep=0}, "{E_I^J = \iota_{I!}^J}", curve={height=-12pt}, hook', from=1-1, to=1-2]
    \arrow[""{name=1, anchor=center, inner sep=0}, "{\Bor_I^J = \iota_I^{J*}}", curve={height=-12pt}, from=1-2, to=1-1]
    \arrow["\dashv"{anchor=center, rotate=-90}, draw=none, from=0, to=1]
  \end{tikzcd}
\]
witnesses $\Op_I \subset \Op_J$ as a colocalizing subcategory.
Moreover, it behaves well with $\obv$.
\begin{proposition}[{\steprop{1.44}}]
  Suppose $\cO^{\otimes}, \cP^{\otimes}$ have at least one color.
  Then, there is an equality 
  \[
    A \prn{\cO \otimes \cP} = A\cO \vee A\cP.
  \]
  In particular, $\Op_I \subset \Op_{\cT}$ is a symmetric monoidal full subcategory.
\end{proposition}

\subsubsection{Restrictions of \tcT-operads}
The $\cT$-category of coefficient systems has a universal property
\[
  \Fun_{\cT}(\cC, \uCoFr^{\cT} \cD) \simeq \Fun(\Tot^{\cT} \cC, \cD)
\]
by \cite[Def.~\href{https://dspace.mit.edu/bitstream/handle/1721.1/112895/1015183829-MIT.pdf?sequence=1\#page=13}{1.10}]{Nardin_thesis} in particular, this yields a \emph{restriction functor}
\[\begin{tikzcd}[ampersand replacement=\&, row sep=tiny]
	{\Fun\prn{\Tot \uSigma_{\cT}, \cS}} \& {\Fun\prn{\Tot \uSigma_{V}, \cS}} \\
	{\Fun_{\cT}\prn{\uSigma_{\cT}, \cS_{\cT}}} \& {\Fun_{V}\prn{\uSigma_{V}, \cS_{V}}}
  \arrow[from=1-1, to=1-2, "\Res_V^{\cT}"]
	\arrow["\simeq"{marking, allow upside down}, draw=none, from=1-1, to=2-1]
	\arrow["\simeq"{marking, allow upside down}, draw=none, from=1-2, to=2-2]
	\arrow[from=2-1, to=2-2]
\end{tikzcd}\]
so that, given a map $W \rightarrow V$ and an $W$-set $S$, we have $\Res_V^{\cT}  \cO(S) \simeq \cO(S)$.
By \EBVsubsec{2.3}, this lifts to a restriction functor on $\cT$-operads
\[\begin{tikzcd}[ampersand replacement=\&, row sep=small, row sep=small, row sep=small, row sep=small, row sep=small, row sep=small, row sep=small, row sep=small, row sep=small]
	{\Op_{\cT}} \& {\Op_V} \\
	{\Fun\prn{\Tot^{\cT}\uSigma_{\cT}, \cS}} \& {\Fun\prn{\Tot^{V}\uSigma_{V}, \cS}}
	\arrow["{\Res_V^{\cT}}", from=1-1, to=1-2]
	\arrow[from=1-1, to=2-1]
	\arrow[from=1-2, to=2-2]
	\arrow[from=2-1, to=2-2]
\end{tikzcd}\]
assembling to an equivalence $\Op_\cT \simeq \Gamma^{\cT} \uOp_{\cT}$;
we will refer to the induced tensor product on $\Op_{\cT}$ as $\obv$.

\subsubsection{$I$-symmetric monoidal categories and $\cO$-algebras}%
\cite{Nardin} constructed a (non-full) subcategory inclusion
\[
  \iota\colon \Cat_I^\otimes \rightarrow \Op_{\cT};
\]
$\cT$-operad maps between $I$-symmetric monoidal categories are called \emph{lax $I$-symmetric monoidal functors}, and morphisms in the image of $\iota$ are called \emph{$I$-symmetric monoidal functors}.

Given $\cO^{\otimes},\cC^{\otimes} \in \Op_{\cT}$, we define \emph{$\cO$-algebras in $\cC^{\otimes}$} to be $\cT$-operad maps $\cO^{\otimes} \rightarrow \cC^{\otimes}$, which naturally fit into an $\infty$-category $\Alg_{\cO}(\cC)$.
These have a \emph{pointwise $\cT$-operad structure} $\uAlg_{\cO}^{\otimes}(\cC)$ given by the internal hom in a presentably symmetric monoidal structure on $\Op_{\cT}$, whose tensor product we write as $\obv$ (see \EBVsubsec{3.2} for the tensor functor and \stesubsec{3.1} for the coherences).
The unit for this symmetric monoidal structure is the $\cT$-operad $\triv^{\otimes}_{\cT} \deq \cN_{I^{\triv} \infty}^{\otimes}$ \EBVprop{3.17}, i.e. there is a canonical equivalence
\begin{equation}\label{eqn:triv}
  \uAlg_{\triv_{\cT}}^{\otimes}(\cO) \simeq \cO^{\otimes}
\end{equation}

Moreover, we verified in \EBVcor{3.13} that whenever $\cC^{\otimes}$ is an $I$-symmetric monoidal $n$-category, $\uAlg_{\cO}^{\otimes}(\cC)$ is as well, and given a $\cT$-operad map $\cO^{\otimes} \rightarrow \cP^{\otimes}$ and an $I$-symmetric monoidal functor $\cC^{\otimes} \rightarrow \cD^{\otimes}$, the induced lax $I$-symmetric monoidal functors
\[
  \uAlg_{\cP}^{\otimes}(\cC) \rightarrow \uAlg_{\cO}^{\otimes}(\cC); \hspace{50pt} \uAlg_{\cO}^{\otimes}(\cC) \rightarrow \uAlg_{\cO}^{\otimes}(\cD) 
\]
are $I$-symmetric monoidal.
In particular, when $\cC^{\otimes}$ is an $I$-symmetric monoidal $\infty$-category and $\cO^{\otimes},\cP^{\otimes}$ are $I$-operads, there are natural $I$-symmetric monoidal equivalence
\begin{equation}\label{eqn:comutativity of alg}
  \uAlg_{\cO}^{\otimes} \uAlg_{\cP}^{\otimes}(\cC) \simeq \uAlg_{\cO \otimes \cP}^{\otimes}(\cC) \simeq \uAlg_{\cP}^{\otimes} \uAlg_{\cO}^{\otimes}(\cC)
\end{equation}
\begin{example}[{\EBVex{3.24}}]\label{ex:e0}
  The $\cT$-operads $\EE_0^{\otimes} \deq \cN_{I_{0,\cT}}^{\otimes}$ and $\EE_\infty^{\otimes} \deq \cN_{I_\infty}^{\otimes}$ are characterized by formulas 
  \[
    \uAlg_{\EE_0}(\cC)_V \simeq \cC_{V, 1_V/}; \hspace{40pt} \uAlg_{\EE_\infty}(\cC)_V \simeq \CAlg(\cC_V).
  \]
  In particular, if $1_V$ is terminal for all $V \in \cT$, then $\uAlg_{\EE_0}(\cC) = \cC_*$ is the $\cT$-category of pointed objects.
\end{example}
The following result will be useful:
we lift the \emph{doctrinal adjunction} to $I$-symmetric monoidal categories.
\begin{proposition}[\stecor{D.5}]\label{Doctrinal adjunction}
  Suppose $L^{\otimes}\colon \cC^{\otimes} \rightarrow \cD^{\otimes}$ is an $I$-symmetric monoidal functor whose underlying $\cT$-functor admits a right adjoint $R$.
  Then, $R$ lifts to a canonical lax $I$-symmetric monoidal right adjoint $R^{\otimes} \vdash L^{\otimes}$.
  Moreover, for any $\cT$-operad $\cO^{\otimes}$ the postcomposition lax $I$-symmetric monoidal functors partake in a lax $I$-symmetric monoidal adjunction
  \[
    L_*^{\otimes}\colon \uAlg_{\cO}^{\otimes}(\cC) \rightleftarrows \uAlg^{\otimes}_{\cO}(\cD)\colon R_*^{\otimes}
  \]
  such that $L^{\otimes}_*$ is $I$-symmetric monoidal.
  If $R^{\otimes}$ is symmetric monoidal then $R_*^{\otimes}$ is symmetric monoidal;
  if $R$ is also fully faithful, then $R_*^{\otimes}$ is fully faithful.
\end{proposition}

\subsubsection{The underlying $\cT$-$\infty$-category}%
Note that composite functor $\Cat_I^{\otimes} \rightarrow \Op_I \xrightarrow{U} \Cat_{\cT}$ is the usual \emph{underlying $\cT$-$\infty$-category} functor.
Moreover, $U$ behaves well with respect to $\uAlg^{\otimes}$;
indeed, we verified in \EBVcor{3.18} that the underlying $\cT$-$\infty$-category has values
\[
  U\prn{\uAlg_{\cO}^{\otimes}(\cC)}_V \simeq \Alg_{\Res_V^{\cT} \cO}\prn{\Res_V^{\cT} \cC}; \hspace{40pt}
  \Alg_{\cO}(\cC) \simeq \Gamma^{\cT} U \uAlg_{\cO}^{\otimes}(\cC).
\]
It was essentially observed in \ns{Cor.}{2.4.5} that the composite functor $\Op_{I^{\triv}} \subset \Op_{\cT} \xrightarrow{U} \Cat_{\cT}$ is an equivalence. 
We write $\triv(-)^{\otimes}$ for the composite functor
\[
  \triv(-)^{\otimes}\colon \Cat_{\cT} \xrightarrow{\;\; U^{-1} \;\;} \Op_{I^{\infty}} \hookrightarrow \Op_{\cT};
\]
\ns{Cor.}{2.4.5} shows that $\triv(\cC)$ algebras are simply $\cC$-indexed diagrams of objects, i.e.
\[
  \uAlg_{\triv(\cC)}(\cO) \simeq \uFun_{\cT}(\cC,U\cO).
\]

\subsubsection{Unital $I$-operads}%
Assume that $I$ is an almost unital weak indexing category.
In \cite{Tensor} we introduced the following gamut of definitions, each of which will be useful.
\begin{definition}
  We say that an $I$-operad $\cO^{\otimes}$ 
  \begin{itemize}
    \item is \emph{almost unital} if it has at least one color and whenever there exists some $S \in \FF_V$ such that $\cO(S) \neq \emptyset$, we have $\cO(\emptyset_V) \simeq *$,
    \item is \emph{unital} if it has at least one color and $\cO(\emptyset_V) \simeq \cN_{I \infty}(\emptyset_V)$ for all $V \in \cT$, and
    \item is \emph{almost reduced} if it is almost unital and has one color, and 
    \item is \emph{reduced} if it is unital and has one color.\qedhere
  \end{itemize}
\end{definition}
A $\cT$-operad is almost unital if and only if it's a unital $I$-operad for \emph{some} almost-unital weak indexing category $I$.
For this reason, we'll usually focus on either unital $I$-operads or almost-unital $\cT$-operads.
It will be important to keep the $I$-symmetric monoidal case in mind.
\begin{example}[\steobs{1.58}]\label{ex:unital I-SMC}
  An $I$-symmetric monoidal $\infty$-category $\cC^{\otimes}$ is a unital $I$-operad if and only if, for all $V \in \upsilon(I)$, the unit object $1_V \in \cC_V$ is initial.
\end{example}

Write $\EE_{0,\upsilon(I)}^{\otimes} \deq \cN_{I_{0,\upsilon(I)}}^{\otimes}$.
We will largely use the following result ti access unital $I$-operads.
\begin{proposition}[{\stecor{2.15}}]\label{thm:unital prop}
The full subcategory $\Op_{I}^{\uni} \subset \Op_I$ of unital $I$-operads is both a localizing and colocalizing subcategory, i.e. the inclusion participates in a double adjunction
\[
  \begin{tikzcd}[ampersand replacement=\&, column sep=large]
    {\Op_I} \& {\Op_I^{\uni}.}
    \arrow[""{name=0, anchor=center, inner sep=0}, "{(-) \obv \EE_{0,\upsilon(I)}^{\otimes}}", curve={height=-18pt}, from=1-1, to=1-2]
    \arrow[""{name=1, anchor=center, inner sep=0}, "{\uAlg_{\EE_{0,\upsilon(I)}}^{\otimes}(-)}"', curve={height=18pt}, from=1-1, to=1-2]
    \arrow[""{name=2, anchor=center, inner sep=0}, hook', from=1-2, to=1-1]
    \arrow["\dashv"{anchor=center, rotate=-90}, draw=none, from=0, to=2]
    \arrow["\dashv"{anchor=center, rotate=-90}, draw=none, from=2, to=1]
  \end{tikzcd}
\]
In particular, if $\cO^{\otimes}$ and $\cC^{\otimes}$ are unital, then there are natural equivalences
\begin{align*}
  \uAlg_{\cP}^{\otimes}(\cC) 
  &\simeq \uAlg_{\cP \otimes \EE_{0,\upsilon(I)}}^{\otimes}(\cC);\\
  \uAlg_{\cO}^{\otimes}(\cD) 
  &\simeq \uAlg^{\otimes}_{\cO} \uAlg^{\otimes}_{\EE_{0,\upsilon(I)}}(\cD).
\end{align*}
\end{proposition}
This and \cref{eqn:triv} together yield the following corollary.
\begin{corollary}\label{thm:E initial}
  $\EE_{0,\upsilon(I)}^{\otimes}$ is initial among reduced $I$-operads.
\end{corollary}

\subsubsection{Algebras in cartesian and cocartesian $I$-symmetric monoidal $\infty$-categories}%
In \steprop{1.51} we gave algebras in cartesian $I$-symmetric monoidal $\infty$-categories an explicit presentation generalizing the \emph{$\cO$-monoids} of \ha{Prop.}{2.4.2.5} (as $\cT$-functors satisfying ``Segal conditions'') which we will not describe explicitly here;
as a relic of this, we will simply use the notation
\begin{equation}\label{Monoid equation 1}
  \uMon_{\cO}(\cD) \deq \uAlg_{\cO}\prn{\cD^{I-\times}}; \hspace{50pt} \Mon_{\cO}(\cD) \deq \Alg_{\cO}\prn{\cD^{I-\times}}.
\end{equation}
The pointwise $I$-symmetric monoidal structure is cartesian \stecor{1.68}.
Given $\cC \in \Cat$, we write
\begin{equation}\label{Monoid equation 2}
  \uMon_{\cO}(\cC) \deq \uMon_{\cO}\prn{\uCoFr^{\cT} \cC}; \hspace{50pt} \Mon_{\cO}(\cC) \deq \Mon_{\cO}\prn{\uCoFr^{\cT} \cC}.
\end{equation}
instead we will care about their monadic presentation, which goes as follows.
\begin{proposition}[\EBVcor{2.73}]\label{Monad prop}
  Suppose $\cC$ is a presentable and cartesian closed $\infty$-category and $\cO^{\otimes}$ a $\cT$-operad with connected colors.
  Then, the monad $T_{\cO}$ associated with the monadic functor $\Mon_{\cO}(\cC) \rightarrow \CoFr^{\cT} \cC$ has fixed points
  \[
    \prn{T_{\cO} X}^W \simeq \coprod_{S \in \FF_{I,W}} \prn{\Fr_{\cC}\cO(S) \times \prod_{U \in \Orb(S)} X^U}_{h\Aut_W(S)},
  \]
  where $\Fr_{\cC}\colon \cS \rightarrow \cC$ is the unique left adjoint sending $*$ to the terminal object of $\cC$.
\end{proposition}

Moreover, in the case that $\cO^{\otimes}$ is unital, we characterized cocartesian algebras simply as diagrams
\begin{equation}\label{Cocartesian algs eq}
  \uAlg_{\cO}^{\otimes}\prn{\cC^{I-\sqcup}} \simeq \uFun_{\cT}(U\cO,\cC)^{I-\sqcup}
\end{equation}
\steprop{1.62};
in fact, $\cC^{I-\sqcup}$ still exists as an $I$-operad with the above algebras in when $\cC$ is not assumed to have $I$-indexed coproducts.
In particular, in the unital case, we acquire a double adjunction
\begin{equation}\label{eq:U adjunction}
  \begin{tikzcd}[ampersand replacement=\&, column sep=large]
    {\Cat_{\cT}} \& {\Op_{I}^{\uni}.}
    \arrow[""{name=0, anchor=center, inner sep=0}, "{\triv(-)^{\otimes} \obv \EE_{0,\upsilon(I)}}", curve={height=-16pt}, from=1-1, to=1-2]
    \arrow[""{name=1, anchor=center, inner sep=0}, "{(-)^{I-\sqcup}}"', curve={height=16pt}, from=1-1, to=1-2]
    \arrow[""{name=2, anchor=center, inner sep=0}, "U"{description}, from=1-2, to=1-1]
    \arrow["\dashv"{anchor=center, rotate=-89}, draw=none, from=0, to=2]
    \arrow["\dashv"{anchor=center, rotate=-91}, draw=none, from=2, to=1]
  \end{tikzcd}
\end{equation}
\begin{example}\label{ex:N infinifty cocart}
  \cref{Cocartesian algs eq} constructs an equivalence
  \[
    \Alg_{\cO}\prn{*_{\cT}^{I-\sqcup}} \simeq * \simeq \Alg_{\cO}\prn{\cN_{I \infty}^{\otimes}},
  \]
  natural in the unital $I$-operad $\cO^{\otimes}$, and hence an equivalence $\cN_{I\infty}^{\otimes} \simeq *_{\cT}^{I-\sqcup}$ via Yoneda's lemma.
\end{example}
\begin{example}[{\stethm{2.2}}]\label{ex:calg-cocart}
  Given $\cC^{\otimes}$ a $\cT$-operad, $\Bor_I^{\cT} \uCAlg_{I}^{\otimes}(\cC)$ is a cocartesian $I$-operad.
\end{example}
This was used in \cite{Tensor} to prove the following important proposition.
\begin{proposition}[{\steprop{2.10}}]\label{prop:tensor-with-cocart}
  Given $\cP^{\otimes}$ an $I$-operad, the canonical map $\cP^{\otimes} \rightarrow \cP^{\otimes} \obv \cN_{I\infty}^{\otimes}$ is an equivalence if and only if there exists an equivalence $\cP^{\otimes} \simeq (U\cP)^{I-\sqcup}$;
  in particular, $\cN_{I\infty}^{\otimes}$ is an idempotent algebra in $\Op_I$.
\end{proposition}

\subsubsection{$I$-$d$-operads}%
In \cite{EBV}, we defined the full subcategory $\Op_{\cT,d} \subset \Op_{\cT}$ of \emph{$\cT$-$d$-operads} to be those such that $\cO(S)$ is a $(d-1)$-truncated space for all $S \in \uFF_{A\cO}$, and verified the following.
\begin{proposition}[{\EBVsubsec{2.5}}]\label{thm:hd}
  Fix $d \geq -1$ and $\cO^{\otimes} \in \Op_{\cT}$.
  \begin{enumerate}
    \item The inclusion $\Op_{\cT,d} \subset \Op_{\cT}$ has a left adjoint $h_d\colon \Op_{\cT} \rightarrow \Op_{\cT, d}$ satisfying
      \[
        h_d \cO(S) \simeq \tau_{\leq d-1} \cO(S).
      \]
    \item The unit of the $h_0$-localization adjunction is the map $\cO^{\otimes} \rightarrow \cN_{A\cO}^{\otimes}$;
      in particular, $\cN_{(-) \infty}^{\otimes}$ factors through an equivalence
      \[
        \wIndCat_{\cT} \simeq \Op_{\cT, 0}.
      \]
    \item When $\cP^{\otimes}$ is a $\cT$-$d$-operad, there is a natural equivalence
      \[
        \uAlg_{\cO}^{\otimes}(\cP) \simeq \uAlg_{h_d\cO}^{\otimes}(\cP),
      \]
      and each are $\cT$-$d$-operads.
    \item An $I$-symmetric monoidal $\infty$-category $\cC^{\otimes}$ is a $\cT$-$d$-operad if and only if $U\cC$ is a $\cT$-$d$-category.
  \end{enumerate}
\end{proposition}
We call $h_d \cO^{\otimes}$ the \emph{homotopy $d$-operad of $\cO^{\otimes}$.}
An important result about this is the following¸ which follows from a multi-colored version of \cref{Monad prop}.
\begin{corollary}[\stecor{A.25}]\label{cor:equivalence on S}
  Suppose $\varphi\colon \cO^{\otimes} \rightarrow \cP^{\otimes}$ is a map of $\cT$-operads inducing an equivalence $\pi_0 U\cO \rightarrow \pi_0 U \cP$.
  Then, $\varphi$ is an $h_{n+1}$-equivalence if and only if the induced functor $\Mon_{\cP}(\cS_{\leq n}) \rightarrow \Mon_{\cO}(\cS_{\leq n})$ is an equivalence.
\end{corollary}

Now, fix $\cC^{\otimes}$ an $I$-symmetric monoidal 1-category;
in light of \cref{thm:hd}, to characterize $\cO$-algebras in $\cC^{\otimes}$, we may equivalently characterise $h_1 \cO$-algebras in $\cC$, so assume $\cO^{\otimes}$ is an $I$-1-operad, i.e. its structure spaces are sets.

We gave a simple combinatorial model for $I$-1-operads in \EBVsubsec{2.7}, which we will not relitigate here, instead focusing only on algebras.
Given a $\cT$-object $X \in \Gamma^{\cT} \cC$, we defined the \emph{unreduced endomorphism $I$-operad} $\End_X(\cC)$ as a one-colored $I$-1-operad with structure sets 
\[
  \End_X(\cC)(S) \simeq \Hom_{\cC_V} \prn{X_{\uV}^{\otimes S}, X_{\uV}},
\]
where $X_{\uV} \in \cC_V$ is the $V$-object underlying $X$.
1-categorical algebras take a familiar form.
\begin{proposition}[{\EBVprop{2.107}}]\label{thm:1-categorical algebras}
  Given $\cO^{\otimes} \in \Op_{I,1}^{\oc}$,
  $\Alg_{\cO}(\cC)$ is a 1-category whose objects are pairs $(X \in \Gamma^{\cT}\cC,\varphi\colon \cO \rightarrow \End_X(\cC))$ and whose morphisms are $\Gamma^{\cT} \cC$-maps $f\colon X \rightarrow Y$ such that the corresponding diagram commutes
\[\begin{tikzcd}[ampersand replacement=\&, row sep=tiny]
	\& {\End_X(\cC)} \\
	{\cO^{\otimes}} \\
	\& {\End_Y(\cC)}
	\arrow["{\End_f}", from=1-2, to=3-2]
	\arrow[from=2-1, to=1-2]
	\arrow[from=2-1, to=3-2]
\end{tikzcd}\]  
\end{proposition}
Moreover, we may exploit this to explicitly describe interchange.
\begin{corollary}[{\EBVobs{3.14}}]\label{thm:1-categorical interchange}
  Given $\cO^{\otimes},\cP^{\otimes} \in \Op_{I,1}^{\oc}$, 
  an $\cO \obv \cP$-algebra structure on $X$ is precisely a pair of $\cO$-algebra and $\cP$-algebra structures such that, for all $\mu \in \cO(S)$, the corresponding $\cC$-map $X^{\otimes S}_{\uV} \rightarrow X_{\uV}$ is a morphism of $\cP$-algebras;
  a morphism of $\cO \obv \cP$-algebras is a $\Gamma^{\cT} \cC$-map which is separately an $\cO$-algebra and $\cP$-algebra morphism. 
\end{corollary}

\section{\tI-operads}
In this section, we establish the $I$-operadic results necessary to prove \cref{LEHA,HEHA,k-connected theorem}.
In particular, in \cref{sec:mapping out} we prove the universal property for algebras \emph{out of} cocartesian $I$-operads, showing compatibility between cocartesian structures and the formation of homotopy $n$-operads.
Using this, in \cref{sec:h equivalences} we give the necessary recognition results on $h_n$-equivalences and $h_n$-cocartesian $I$-operads for the rest of the paper.
Last, in \cref{sec:reduced} we characterize the reduced endomorphism $I$-operad, ultimately forming the main content of the reduction of \cref{LEHA,HEHA} to \cref{k-connected theorem}.
\subsection{The mapping-out property for cocartesian structures}\label{sec:mapping out}
\begin{proposition}\label{Algebras over cocartesian thing}
  Given $\cC \in \Cat_{\cT}$, $\cC^{I-\sqcup}$ is identified by the mapping-out property 
  \[
    \Alg_{\cC^{I-\sqcup}}(\cD) \simeq \Fun_{\cT}\prn{\cC, \uCAlg_I(\cD)};
  \]
  in particular, we acquire a triple adjunction
  \begin{equation}\label{eq:triple adjunction}
    \begin{tikzcd}[ampersand replacement=\&]
	{\Cat_{\cT}} \&\& {\Op_{I}^{\uni}.}
  \arrow[""{name=0, anchor=center, inner sep=0}, "{{\triv(-)^{\otimes} \obv \EE_{0,\upsilon(I)}^{\otimes}}}", curve={height=-30pt}, from=1-1, to=1-3, hook']
	\arrow[hook, ""{name=1, anchor=center, inner sep=0}, "{{(-)^{I-\sqcup}}}"{description}, curve={height=10pt}, from=1-1, to=1-3]
	\arrow[""{name=2, anchor=center, inner sep=0}, "U"{description}, curve={height=10pt}, from=1-3, to=1-1]
	\arrow[""{name=3, anchor=center, inner sep=0}, "{\uCAlg_I^{\otimes}(-)}", curve={height=-30pt}, from=1-3, to=1-1]
	\arrow["\dashv"{anchor=center, rotate=-89}, draw=none, from=1, to=3]
	\arrow["\dashv"{anchor=center, rotate=-91}, draw=none, from=0, to=2]
	\arrow["\dashv"{anchor=center, rotate=-90}, draw=none, from=2, to=1]
\end{tikzcd}
\end{equation}
such that, $\triv(-)^{\otimes} \obv \EE_{0,\upsilon(I)}$ and $(-)^{I-\sqcup}$ are fully faithful.
\end{proposition}
\begin{proof}
  For the first statement, simply apply the equivalences
  \begin{align*}
    \Alg_{\cC^{I-\sqcup}}(\cD) 
    &\simeq \Alg_{\cC^{I-\sqcup} \otimes \cN_{I\infty}}(\cD)
    & \text{\cref{prop:tensor-with-cocart}},\\
    &\simeq \Alg_{\cC^{I-\sqcup}} \uCAlg_I^{I-\sqcup}(\cD) & \text{\cref{ex:calg-cocart,eqn:comutativity of alg}}\\
    &\simeq \Fun_{\cT}\prn{\cC,\uCAlg_I(\cD)} 
    & \text{\cref{Cocartesian algs eq}.}    
  \end{align*}
and Yoneda's lemma under the equivalence $\Alg_{\cO}(\cP)^{\simeq} \simeq \Map_{\Op_{\cT}}\prn{\cO^{\otimes}, \cP^{\otimes}}$.
The bottom adjunction follows by taking cores, and the remaining adjunctions by \cref{eq:U adjunction}.
Fully faithfulness for the former follows from the latter, which is itself follows by combining the mapping out property with \cref{Cocartesian algs eq} and taking cores.
\end{proof}

\begin{remark}
  The case $I = \cT$ is proved in \cite[Lem.~\href{https://arxiv.org/pdf/2503.03024v1\#equation.4.1.10}{4.1.10}]{Yang_HKR}, but it is used as \emph{input to} rather than a corollary of computations of cocartesian algebras, so their techniques are more difficult.
\end{remark}
  \begin{remark}
    \cref{Algebras over cocartesian thing} is the restriction of the \emph{defining} property of $\uCAlg_{\cT}^I(\cC)$ in \cite{Lenz_norms}, who left implicit a comparison with the atomic orbital setting;
  \cref{Algebras over cocartesian thing} together with the identification of the two notions of cocartesian structure gives a slick identification of the two \cite{Tensor}.
\end{remark}

We easily acquire compatibility of $h_n$ with cocartesian structures.
\begin{corollary}\label{h_n commutes lemma}
  Given $\cC$ a $\cT$-category, there exists an equivalence $h_n \prn{\cC^{I-\sqcup}} \simeq \prn{h_n \cC}^{I-\sqcup}$.
\end{corollary}
\begin{proof}
  \cref{thm:hd} constructs a commutative diagram
  \[\begin{tikzcd}[ampersand replacement=\&]
	{\Op_{\cT, d}} \& {\Cat_{\cT,d}} \\
	{\Op_{\cT}} \& {\Cat_{\cT}}
	\arrow["{\uCAlg_{I}(-)}", from=1-1, to=1-2]
	\arrow[hook, from=1-1, to=2-1]
	\arrow[hook', from=1-2, to=2-2]
	\arrow["{\uCAlg_{I}(-)}", from=2-1, to=2-2]
\end{tikzcd}\]
  The result follows by taking left adjoints.
\end{proof}

\subsection{Recognizing \tI-local \texorpdfstring{$h_{n}$}{h}-equivalences}\label{sec:h equivalences}
\subsubsection{Detection via algebras}
\cref{k-connected theorem} will recognize some morphisms of $\cT$-operads which become equivalences after applying $h_{n+1}$, so we now spell out some of its antecedents.
\begin{proposition}\label{Mapping spaces prop}
  Let $\varphi\cln \cO^{\otimes} \rightarrow \cP^{\otimes}$ be a morphism of $\cT$-operads inducing an equivalence $\pi_0 U \cO \rightarrow \pi_0 U \cP$.
  The following are equivalent:
  \begin{enumerate}[label={(\alph*)}]
    \item \label[condition]{cond:equivalence} $\varphi$ is an $h_{n+1}$-equivalence;
    \item \label[condition]{cond:equivalence on C} for all $\cT$-symmetric monoidal $(n+1)$-categories $\cC$, the pullback $\cT$-symmetric monoidal functor
      \[
        \uAlg^{\otimes}_{\cP}(\cC) \rightarrow \uAlg_{\cO}^{\otimes}(\cC)
      \]
      is an equivalence;
    \item \label[condition]{cond:equivalence on S} the pullback functor
      \[
        \Mon_{\cP}(\cS_{\leq n}) \rightarrow \Mon_{\cO}(\cS_{\leq n})
      \]
      is an equivalence; and
    \item \label[condition]{cond:equivalence on presheaves} for all $\infty$-categories $K$, the pullback map of spaces
      \[
        \Mon_{\cP}\prn{\cS_{\leq n}^K}^{\simeq}
        \rightarrow \Mon_{\cO}\prn{\cS_{\leq n}^K}^{\simeq}
      \]
      is an equivalence.
  \end{enumerate}
\end{proposition}
Analogously to \SY{Prop.}{3.2.6}{22}, to prove this, we will apply the following lemma.
\begin{lemma}\label{Monoids in presheaves}
  Given a $\cT$-operad $\cP^{\otimes}$ and a pair of $\infty$-categories $\cD,K$ such that $\cD$ admits finite products, there is an equivalence
  \[
    \uMon_{\cP} \prn{\cD^K} \simeq \uFun_{\cT} \prn{\Infl_e^{\cT} K, \uMon_{\cP}(\cD)}, 
  \]
  natural in functors of $K$, product-preserving functors of $\cD$, and $\cT$-operad maps of $\cP$;
  in particular, taking $\cT$-fixed points yields a natural equivalence of categories
  \[
    \Mon_{\cP} \prn{\cD^K} \simeq \Mon_{\cP}(\cD)^K.
  \]
\end{lemma}
\begin{proof}
  We construct a chain of equivalences
  \begin{align*}
    \uMon_{\cP}\prn{\cD^K}
    &\simeq \uAlg_{\cP}(\uCoFr^{\cT}(\cD^K)^{\cT-\times})  
      &\text{\cref{Monoid equation 1,Monoid equation 2}}\\
    &\simeq \uAlg_{\cP} \uFun_{\cT}\prn{\Infl_e^{\cT} K, \uCoFr^{\cT} \cD}^{\cT-\times}   
      &\text{\cref{Inflation example}}\\
    &\simeq \uAlg_{\cP} \uAlg^{\otimes}_{\triv(\Infl_{e}^{\cT} K)}\prn{\uCoFr^{\cT} \cD^{\cT-\times}} 
      &\text{\cref{eqn:triv}}\\
      &\simeq \uAlg_{\triv(\Infl_e^{\cT} K)} \uAlg_{\cP}^{\otimes} \prn{\uCoFr^{\cT} \cD^{\cT-\times}} 
      &\text{\cref{eqn:comutativity of alg}}\\
    &\simeq \uFun_{\cT} \prn{\Infl_e^{\cT} K, \uAlg_{\cP} \prn{\uCoFr^{\cT}, \cD^{\cT-\times}}} 
      &\text{\cref{eqn:triv}}\\
    &\simeq  \uFun_{\cT} \prn{\Infl_e^{\cT} K, \uMon_{\cP}\prn{\cD}} 
    &\text{\cref{Monoid equation 1,Monoid equation 2}}
  \end{align*}
  The remaining equivalence follows by noting that $\Gamma^{\cT} \Infl_e^{\cT} \cC \simeq \cC$, naturally in $\cC$.
\end{proof}

\begin{proof}[Proof of \cref{Mapping spaces prop}]
  The implication \cref{cond:equivalence} $\implies$ \cref{cond:equivalence on C} is \cref{thm:hd}.
  Moreover, the implications \cref{cond:equivalence on C} $\implies $ \cref{cond:equivalence on S,cond:equivalence on presheaves} is obvious.
  The implication \cref{cond:equivalence on S} $\implies$ \cref{cond:equivalence} is \cref{cor:equivalence on S}.
  Moreover, fixing $\cD = \cS_{\leq n}$ and taking cores of \cref{Monoids in presheaves} yields a natural equivalence
  \[
    \Mon_{\cP}\prn{\cS_{\leq n}^K}^{\simeq} \simeq \Map_{\Cat}\prn{K, \Mon_{\cP}\prn{\cS_{\leq n}}}
  \]
  so \cref{cond:equivalence on presheaves} and Yoneda's lemma together imply \cref{cond:equivalence on S}.
\end{proof}

\subsubsection{The smashing localization on $\cT$-$n$-operads associated with $\cN_{I\infty}^{\otimes}$}
Note the following.
\begin{proposition}\label{hn compatible prop}
  If $\varphi\colon \cO^{\otimes} \rightarrow \cP^{\otimes}$ is an $h_n$-equivalence and $\cQ^{\otimes}$ is a $\cT$-operad, then the induced map 
  \[
    \cQ^{\otimes} \obv \varphi\colon \cQ^{\otimes} \obv \cO^{\otimes} \longrightarrow \cQ^{\otimes} \obv \cP^{\otimes}
  \]
  is an $h_n$-equivalence.
\end{proposition}
\begin{proof}
  By \cref{Mapping spaces prop}, pullback along $\varphi \otimes \cQ^{\otimes}$ yields an equivalence
  \[
    \begin{tikzcd}[ampersand replacement=\&, row sep=tiny]
      {\Mon_{\cQ} \uMon_{\cP}\prn{\cS_{\leq n}}} \& {\Mon_{\cQ} \uMon_{\cO}\prn{\cS_{\leq n}}} \\
      {\Mon_{\cQ \otimes \cP}\prn{\cS_{\leq n}}} \& {\Mon_{\cQ \otimes \cO}\prn{\cS_{\leq n}}}
      \arrow[from=1-1, to=1-2]
      \arrow["\simeq"{marking, allow upside down}, draw=none, from=1-1, to=2-1]
      \arrow["\simeq"{marking, allow upside down}, draw=none, from=1-2, to=2-2]
      \arrow[from=2-1, to=2-2]
    \end{tikzcd}
  \]
  Noting that $\pi_0 U \varphi$ is an equivalence, applying \cref{Mapping spaces prop} shows that $\varphi \otimes \cQ^{\otimes}$ is an $h_n$-equivalence.
\end{proof}
In particular, \cref{hn compatible prop} and \ha{Prop.}{2.2.1.8} construct a symmetric monoidal structure on $\Op_{\cT, n}$ together with a symmetric monoidal structure on $h_{n}$.
The tensor product for this structure is $\cO^{\otimes} \obv_n \cP^{\otimes} \simeq h_n \prn{\cO^{\otimes} \obv \cP^{\otimes}}$, and in particular, \cref{prop:tensor-with-cocart} shows that $\cN_{I \infty}^{\otimes} \in \Op_{\cT, n}$ is an idempotent algebra.
It's easy to identify its smashing localization, and in fact, its $h_n$-preimages.
\begin{corollary}\label{hn smashing corollary}
  Suppose $\cO^{\otimes}$ is an almost-unital $\cT$-operad.
  Then, the following conditions are equivalent:
  \begin{enumerate}[label={(\alph*')}]
    \item The map $\Bor_I^{\cT} \cO^{\otimes} \rightarrow \prn{h_n U\cO}^{I-\sqcup}$ is an $h_n$-equivalence. 
      \label[condition]{cond:n connected at I} 
      \setcounter{enumi}{4}
    \item \label[condition]{cond:op cocart} For all $\cT$-$(n+1)$-operads $\cP^{\otimes}$, the $\cT$-operad $\uAlg_{\cO}^{\otimes}(\cP)$ is cocartesian at $I$.
    \item \label[condition]{cond:cocart} For all $A\cO$-symmetric monoidal $(n+1)$-categories $\cC^{\otimes}$, the $A\cO$-symmetric monoidal $(n+1)$-category $\uAlg_{\cO}^{\otimes}(\cC)$ is cocartesian at $I$.
    \item \label[condition]{cond:semi} The $\cT$-$(n+1)$-category $\uMon_{\cO}\prn{\cS_{\leq n}}$ is $I$-semiadditive.
      \setcounter{enumi}{9}
    \item \label[condition]{cond:smash} The unit map tensors to an $h_n$-equivalence
      \[
        h_n \prn{\id \otimes !}\colon h_n \cO^{\otimes} \simeq h_n \prn{\cO^{\otimes} \obv \triv_{\cT}^{\otimes}} \xrightarrow{\;\;\;\;\; \sim \;\;\;\;\;} h_n\prn{\cO^{\otimes} \obv \cN_{I \infty}^{\otimes}}.
      \]
  \end{enumerate}
\end{corollary}
\begin{proof}
  The implication \cref{cond:n connected at I} $\implies$ \cref{cond:op cocart} is \cref{thm:hd};
  the implication \cref{cond:op cocart} $\implies$ \cref{cond:cocart} is obvious;
  the implication \cref{cond:cocart} $\implies$ \cref{cond:semi} is \cref{cor:cart} and cartesianness of the pointwise structure on $\uMon_{\cO}(\cC)$.
  The implication \cref{cond:semi} $\implies$ \cref{cond:smash} follows straightforwardly from \cref{cor:equivalence on S} and \cref{ex:calg-cocart};
  this was also spelled out in \stecor{2.4}.
  Lastly, the implication 
  \cref{cond:smash} $\implies$ \cref{cond:n connected at I} follows from the following:
  \begin{align*}
    h_n \Bor_I^{\cT}  \prn{\cO^{\otimes}} 
    &\simeq h_n \Bor_I^{\cT} \prn{\cO^{\otimes}  \obv \cN_{I\infty}} & \text{\cref{cond:smash}}\\
    &\simeq h_n \prn{U\cO}^{I-\sqcup} & \text{\cref{prop:tensor-with-cocart}}\\
    &\simeq \prn{h_nU\cO}^{I-\sqcup} & \text{\cref{h_n commutes lemma}}.&\qedhere
  \end{align*}
\end{proof}

\subsection{The reduced endomorphism \tI-operad as a right adjoint}\label{sec:reduced}
In \cite{Tensor}, we introduced the \emph{reduced endomorphism $I$-operad} of a $\cT$-operad for the purpose of lifting the disintegration and assembly process of \cite{HA}.
In this section, we gain explicit computational control over reduced endomorphism $I$-operads of unital $I$-symmetric monoidal $\infty$-categories as follows.
\begin{proposition}\label{thm:reduced endomorphism I-operad}
  The inclusion $\Op_{I}^{\red} \simeq \Op_{I, \EE_{0,\upsilon(I)}/}^{\red} \hookrightarrow \Op_{I, \EE_{0,\upsilon(I)}/}^{\uni}$ has a right adjoint computed by the pullback
  \begin{equation}\label{Reduced endomorphism I-operad equation}
    \begin{tikzcd}[ampersand replacement=\&]
      {\End^{I,\red}_X} \& {\cO^{\otimes}} \\
      {\cN_{I \infty}^{\otimes}} \& {\cO^{\cT-\sqcup}}
      \arrow[from=1-1, to=1-2]
      \arrow[from=1-1, to=2-1]
      \arrow["\lrcorner"{anchor=center, pos=0.125}, draw=none, from=1-1, to=2-2]
      \arrow["\eta", from=1-2, to=2-2]
      \arrow["{\cbr{X}}"', from=2-1, to=2-2]
    \end{tikzcd}
  \end{equation}
  In the case that $\cC^{\otimes}$ is a unital $I$-symmetric monoidal $\infty$-category and $X \in \cC_V$ is a $V$-object,
  mapping in from the free unital $I$-operad $\Fr_{\Op}(S) \obv \EE_{0,\upsilon(I)}$ on an operation in arity $S \in \FF_{I,V}$ yields a pullback
  \[
    \begin{tikzcd}[ampersand replacement=\&]
      {\End^{I,\red}_X(S)} \& {\Map_{\cC_V}\prn{X^{\otimes  S}, X}} \\
      {\cbr{\nabla}} \& {\Map_{\cC_V}\prn{X^{\sqcup  S}, X}}
      \arrow[from=1-1, to=1-2]
      \arrow[from=1-1, to=2-1]
      \arrow["\lrcorner"{anchor=center, pos=0.125}, draw=none, from=1-1, to=2-2]
      \arrow[from=1-2, to=2-2, "W_{S,X}^*"]
      \arrow[from=2-1, to=2-2]
    \end{tikzcd}
  \]
  i.e. 
  $\End^{I,\red}_X(S)$ is equivalent to the space of lifts along the following dashed arrow in $\cC_V$
  \[
    \begin{tikzcd}[ampersand replacement=\&]
      {X^{\sqcup S}} \& X \\
      {X^{\otimes S}} \& {*}
      \arrow["\nabla", from=1-1, to=1-2]
      \arrow["W_{S,X}"', from=1-1, to=2-1]
      \arrow["{!}", from=1-2, to=2-2]
      \arrow[dashed, from=2-1, to=1-2]
      \arrow["{!}"', from=2-1, to=2-2]
    \end{tikzcd}
  \]
\end{proposition}
\begin{proof}
  We will apply the general reduction procedure of \SY{Prop.}{2.1.5}{9} to the \emph{sliced} adjunction 
  \[
    U_*\colon \Op_{I, \EE_{0,\upsilon(I)}^{\otimes}/}^{\uni} \longrightleftarrows \Cat_{\cT,*}\colon \eta^* (-)^{I-\sqcup},
  \]
  whose right adjoint is $(-)^{I-\sqcup}$ together with the precomposed structure map
  \[
    \EE_{0,\upsilon(I)}^{\otimes} \xrightarrow{\eta} \cN_{I\infty}^{\otimes} \simeq *_{\cT}^{I-\sqcup} \rightarrow \cC^{I-\sqcup}.
  \] 
  Indeed, $\Cat_{\cT,*}$ admits an initial object $*_{\cT} \simeq U\EE_{0,\upsilon(I)}$ and $\Op_{I,\EE^{\otimes}_{0,\upsilon(I)}/}^{\otimes}$ admits all limits, which are preserved by $U$ since it is a right adjoint by \cref{eq:U adjunction}.
  Moreover, $\EE_{0,\upsilon(I)} \in \Op_{I}^{\red}$ is initial by \cref{thm:E initial}, there is a unique equivalence $\cN_{I \infty}^{\otimes} \simeq *_{\cT}^{I-\sqcup}$ by \cref{ex:N infinifty cocart}, and $\cO^{\otimes} \in \Op_{I, \EE_{0,\upsilon(I)}/}^{\uni}$ corresponds with a reduced $I$-operad if and only if $U\cO^{\otimes} \in \Cat_{\cT,*}$ is initial, so the first claim follows by \SY{Prop.}{2.1.5}{9}.
 
  To acquire the second pullback square, one need only note that the natural equivalences
  \begin{align*}
    \Map_{\Op_{\cT}}\prn{\Fr_{\Op}(S) \obv \EE_{0,\upsilon(I)}, \; \cC^{\otimes}} &\simeq \Map_{\cC_V}\prn{X^{\otimes  S}, X}, \\
    \Map_{\Op_{\cT}}\prn{\Fr_{\Op}(S) \obv \EE_{0, \upsilon(I)}, \; \cN_{I \infty}^{\otimes}} &\simeq *
  \end{align*}
  follow by \cref{thm:free operads,thm:unital prop}.
  What remains is to verify that the $\eta$ induces $W_{S,X}^*$ and the bottom arrow includes the fold map $\nabla$;
  both facts were verified in \stesubsec{A.5}.
\end{proof}

In fact, \SY{Prop.}{4.2.8}{30} introduced a result on connectivity of such spaces of lifts, immediately yielding the following corollary.
\begin{corollary}\label{Schlank connectivity prop}
  If $X \in \cC_V$ is a $(k + \ell + 2)$-truncated object and the Wirthm\"uller map $W_{S,X}\colon X^{\sqcup S} \rightarrow X^{\otimes S}$ is $\ell$-connected, then the space $\End_X^{I,\red}(\cC)(S)$ is $k$-truncated.
\end{corollary}

Given a reduced $I$-operad $\cP^{\otimes}$ and a $V$-object $X \in \cO_V$, applying $\cP$-algebras to \cref{Reduced endomorphism I-operad equation} yields a pullback 
\begin{equation}\label{Disintegration fiber square}
  \begin{tikzcd}[ampersand replacement=\&]
    {\Alg_{\Res_V^{\cT} \cP} \End_X^{I,\red}(\cO)} \& {\uAlg_{\cP}(\cO)_V} \\
    {\cbr{X}} \& {U\cO_V}
    \arrow[from=1-1, to=1-2]
    \arrow[from=1-1, to=2-1]
    \arrow["\lrcorner"{anchor=center, pos=0.125}, draw=none, from=1-1, to=2-2]
    \arrow[from=1-2, to=2-2, "U"]
    \arrow[hook, from=2-1, to=2-2]
  \end{tikzcd}
\end{equation}
In the case that $U\cO$ is a $\cT$-space, $U$ is a automatically cocartesian fibration.
Unfortunately, this is far from our case;
the best we can do is take cores of the above pullback square, resulting in the following proposition.
\begin{proposition}\label{Checking morita equivalences on End}
  Suppose $\cP^{\otimes} \rightarrow \cQ^{\otimes}$ is a morphism of $I$-operads inducing an equivalence of spaces
  \[
    \varphi^{*,\simeq}_{X}\colon \Alg_{\Res_V^{\cT} \cQ} \End_X^{I,\red}(\cO)^{\simeq}
    \longrightarrow \Alg_{\Res_V^{\cT} \cP} \End_X^{I,\red}(\cO)^{\simeq}
  \]
  for all $V \in \cT$ and $X \in U\cO_V$.
  Then, the induced map of $\cT$-spaces
  \[
    \uAlg_{\cQ}(\cO)^{\simeq} \rightarrow \uAlg_{\cP}(\cO)^{\simeq}
  \]
  is an equivalence;
  in particular, passing to $\cT$-fixed points, the induced map of spaces
  \[
    \Alg_{\cQ}(\cO)^{\simeq} \rightarrow \Alg_{\cP}(\cO)^{\simeq}
  \]
  is an equivalence.
\end{proposition}
\begin{proof}
  Taking cores of \cref{Disintegration fiber square}, we find that that $\varphi^{*, \simeq}_{X}$ is the induced map on the homotopy fiber over $X$ of the following map of $\cT$-spaces over $U\cO$:
  \[\begin{tikzcd}[ampersand replacement=\&, row sep=small]
	{\uAlg_{\cQ}(\cO)^{\simeq}} \&\& {\uAlg_{\cP}(\cO)^{\simeq}} \\
	\& {U\cO}
	\arrow["{\varphi^{*,\simeq}}", from=1-1, to=1-3]
	\arrow[from=1-1, to=2-2]
	\arrow[from=1-3, to=2-2]
\end{tikzcd}\]
  $\varphi^{*,\simeq}$ is an equivalence if and only if its $V$-fixed points are an equivalence for all $V \in \cT$, and the homotopy fibers of $\varphi^{*,\simeq,V}$ are contractible by the above argument, so $\varphi^{*, \simeq, V}$ is an equivalence for all $V$. 
  Hence $\varphi^{*, \simeq}$ is an equivalence, proving the proposition.
\end{proof}

\section{Connectivity and Eckmann-Hilton arguments}
We now prove \cref{LEHA,HEHA,k-connected theorem}, beginning with a recognition result for $\ell$-connected $\cO$-monoid maps.
\subsection{Connectivity of algebras can be detected in the value topos}
Fix $\cC$ an $n$-topos for some $n \leq \infty$.
\begin{lemma}\label{Elementwise connectivity}
  A map $f\colon C \rightarrow D$ in $\CoFr^{\cT} \cC$ is $\ell$-connected if and only if, for all $V \in \cT^{\op}$, the fixed point map $C^V \rightarrow D^V$ is $\ell$-connected.
\end{lemma}
\begin{proof}
  Per \cref{Topos-theoretic connectedness}, it is equivalent to prove that $\ell$-connectiveness of a morphism in $\Fun(\cT^{\op}, \cC)$ is measured elementwise.
  Indeed, since (co)limits in $\Fun(\cT^{\op}, \cC)$ are computed elementwise, effective epimorphisms and diagonals are as well.
  The former proves the statement for $(-2)$-connectiveness, and the latter together with the diagonal presentation of \httprop{6.5.1.18} 
  shows that the statement for $(\ell - 1)$-connectiveness implies the statement for $\ell$-connectiveness, so the lemma follows by induction.
\end{proof}

\begin{proposition}\label{U conservativity}
  Given a map $f\colon X \rightarrow Y$ in $\Mon_{\cO}(\cC)$, if the underlying map $Uf$ in $\CoFr^{\cT} \cC$ is $\ell$-connected, then $f$ is $\ell$-connected.
\end{proposition}
\begin{proof}
  In view of \SY{Lem.}{4.4.1}{33} and \cref{Monad prop}, it suffices to verify that the monad $T_{\cO}\colon \CoFr^{\cT} \cC \rightarrow \CoFr^{\cT} \cC$ preserves $\ell$-connected morphisms;
  by \cref{Elementwise connectivity}, it suffices to verify that whenever each $\cC$-diagram $X^V \rightarrow Y^V$ is $\ell$-connected, each induced map $T_{\cO} X^W \rightarrow T_{\cO} X^W$ is $\ell$-connected.
  But by the monad formula of \cref{Monad prop}, it suffices to note that $\ell$-connected morphisms in an $\infty$-topos are closed under cartesian products and colimits \httcorprop{6.5.1.13}{5.2.8.6}.
\end{proof}
For instance, $U$ preserves the terminal object and is conservative, so it also reflects the property of being terminal;
applying \cref{U conservativity} in the case $Y = *$ shows that $U$ reflects $n$-connectivity of objects.
\begin{remark}
  Since $U$ is a right adjoint, it preserves $n$-truncatedness and $n$-truncated objects.
\end{remark}
\begin{warning}
  \cref{U conservativity} is delicate for a few reasons.
  \begin{enumerate}
    \item 1f $\cO$ is not $n$-connected, then the free $\cO$-algebra monad $T_{\cO}\colon \cC_V \rightarrow \cC_V$ may itself may fail to preserve $n$-connected objects;
      indeed, we have $T_{\cO} *_V \simeq \coprod_{S \in \FF_V} \Fr_{\cC} \cO(S)_{h\Aut_V S}$, which is often not much more highly connected than the individual spaces $\cO(S)_{h\Aut_V S}$.
    \item $U$ does not generally \emph{preserve} $\ell$-connectivity of objects or morphisms
      for instance, given an $\ell \geq (k+1)$-connected space $X$, the equivalence $\Omega^k\colon \cS_{*, \geq k+1} \xrightarrow{\; \sim \;} \Alg_{\EE_k}(\cS_{\geq 1})$ of \cite{Guillou-May,Juran} exhibits $\Omega^kX$ as an $\ell$-connected $\EE_k$-algebra such that $U\Omega^nX$ is only in general $(\ell - k)$-connected.
    \item For a similar reason, $U$ does not usually reflect $\ell$-truncatedness of morphisms or objects.\qedhere 
  \end{enumerate}
\end{warning}

\subsection{The proof of \texorpdfstring{\cref{k-connected theorem}}{Theorem D}}%
We now begin to reduce \cref{k-connected theorem} to the case $n \leq \ell + 1$ with the following.
\begin{lemma}\label{Wirthmuller truncation lemma}
  The truncation functor $\tau_{\leq \ell}\colon \cC \rightarrow \tau_{\leq \ell} \cC$ extends to a $\cT$-functor
  \[
    \tau_{\cO}\colon \uMon_{\cO}(\cC) \rightarrow \uMon_{\cO}(\tau_{\leq \ell} \cC)
  \]
  satisfying $\tau_{\cO} W_{S,X}=W_{S,\tau_{\cO} X}$.
  Moreover, the inclusion $\iota\colon \tau_{\leq \ell}\cC \rightarrow \cC$ extends to a fully faithful $\cT$-functor
  \[
    \iota_{\cO}\colon \uMon_{\cO}(\tau_{\leq \ell} \cC) \hookrightarrow \uMon_{\cO}(\cC)
  \]
  such that $\tau_{\cO} W_{S, \iota_{\cO}X} = W_{S,X}$.
\end{lemma}
\begin{proof}
  Since $\tau_{\leq \ell}$ is product-preserving \httlem{6.5.1.2},
  $\tau_{\leq \ell}\colon \uCoFr^{\cT} \cC \rightarrow \uCoFr^{\cT} \tau_{\leq \ell \cC}$ is a $\cT$-symmetric monoidal left adjoint for the cartesian structure \cite{Tensor}; 
  everything other than the equalities involving $W_{S,X}$ then follows straightforwardly from \cref{Doctrinal adjunction}.

  In particular, $\tau_{\cO}$ is a $\cT$-functor which preserves indexed products and coproducts;
  this implies that $\tau_{\cO} W_{S,X} = W_{S,\tau_{\cO} X}$.
  The remaining equality follows from fully faithfulness by noting that
  \[
    \tau_{\cO} W_{S, \iota_{\cO} X} = W_{S, \tau_{\cO} \iota_{\cO} X} = W_{S,X}.\qedhere
  \]
\end{proof}
We say that a map $f\colon X \rightarrow Y$ in an $n$-topos is an \emph{$\ell$-equivalence} if it is a $\tau_{\leq \ell}$-equivalence;
if $f$ admits a section, this is equivalent to $f$ being $\ell$-connected (see \SY{Prop.}{4.3.5}{32} or note that this follows by splitting the long exact sequence in homotopy).
In general, it is \emph{implied} by $\ell$-connectedness, as $\tau_{\leq \ell}$ preserves $\ell$-connectedness and $\ell$-connected maps between $\ell$-truncated objects are equivalences.
We apply this by equivariantizing \SY{Lem.}{5.1.1}{35}.

\begin{lemma}\label{section lemma}
  If $\cC^{I-\times}$ is a Cartesian $I$-symmetric monoidal $\infty$-category and $S \in \uFF_{I}$, then the image of the $\cO$-algebra Wirthm\"uller map $W_{X,S}\colon \coprod^S_U X_U \rightarrow \prod^S_U X_U$ under $U\colon \Alg_{\cO}(\cC)_V \rightarrow \cC_V$ admits a section.
\end{lemma}

\begin{proof}[Proof of \cref{section lemma}]
  Let $i_U\colon Y_U \rightarrow \Res_U^V \coprod_{U'}^S Y_{U'}$ be adjunct to the inclusion $\Ind_U^V Y_U \hookrightarrow \coprod_{U'}^S Y_{U'}$ and fix an operation $\mu \in \cO(S)$.
  We verify commutativity of the following diagram, giving a section $\mu_{\coprod X_U} \sigma_1 f$ for $W_{X,S}$.
    \[\begin{tikzcd}[ampersand replacement=\&, column sep=8em]
    {\prod\limits_U^S \prn{\Res_U^V\coprod\limits_U^V X_U}} \& {\prn{\coprod\limits_U^S X_U}^{\times S}} \& {\coprod\limits_U^S X_U} \\
	\& {\prn{\prod\limits_U^S X_U}^{\times S}} \& {\prod\limits_U^S X_U} \\
	{\prod\limits_U^S X_U} \& {\prod\limits_U^S X_U^{\times \Res_U^V S}} \& {\prod\limits_U^S X_U}
	\arrow["\sim"', "\sigma_1", from=1-1, to=1-2]
  \arrow["{{{\mu_{\coprod X_U}}}}", from=1-2, to=1-3]
	\arrow["{{{h = \prn{W_{\Res_U^V X, \Res_U^V S}}_{U \in \Orb(S)}}}}", from=1-2, to=2-2]
	\arrow["{{W_{X,S}}}", from=1-3, to=2-3]
  \arrow["{{{\mu_{\prod X_U}}}}", from=2-2, to=2-3]
	\arrow["{\sigma_2}"', "\sim", from=2-2, to=3-2]
	\arrow[equals, from=2-3, to=3-3]
  \arrow["{{f = (i_U)_{U \in \Orb(S)}}}", from=3-1, to=1-1]
  \arrow["{{{g = (\iota_U)_{U \in \Orb(S)}}}}", from=3-1, to=3-2]
	\arrow[curve={height=24pt}, equals, from=3-1, to=3-3]
  \arrow["{\prn{\Res_U^V \mu_{X_U}}_{U \in \Orb(S)}}", from=3-2, to=3-3]
    \end{tikzcd}\]
    Here, $\mu_X\colon X^S \rightarrow X$ is the structure map corresponding with $\mu$ for a $\Res^{\cT}_V \cO$-algebra $X$.
    
    Note that the top right square is commutative by the fact that $W_{S,X}$ is an $\cO$-algebra morphism and the bottom right follows by unwinding definitions.  
    Moreover, $\mu \circ g$ is the external product of a collection of endomorphisms $X_U \xrightarrow{\iota_U} X_U^{\times \Res_U^V S} \xrightarrow{\mu} X_U$, the first of which is 
    the inclusion of a unit on all but one factor:
    \[
      \begin{tikzcd}[ampersand replacement=\&, column sep=huge, row sep=small]
	{X_U} \& {X_U^{\times \Res_U^V S}} \& {X_U} \\
	{X_U \times \prod\limits_{W}^{\Res_U^V S - \cbr{\alpha}} 1_W} \& {X_U \times \prod\limits_{W}^{\Res_U^V S - \cbr{\alpha}} X_W}
	\arrow["{\iota_U}", from=1-1, to=1-2]
	\arrow["\simeq"{marking, allow upside down}, draw=none, from=1-1, to=2-1]
  \arrow["\Res_U^V \mu", from=1-2, to=1-3]
	\arrow["\simeq"{marking, allow upside down}, draw=none, from=1-2, to=2-2]
	\arrow["{{(\id,\eta)}}", from=2-1, to=2-2]
	\arrow[from=2-2, to=1-3]
\end{tikzcd}
    \]
    in particular, $\mu \circ \iota_U$ is homotopic to the identity by unitality, so $\mu \circ g \sim \id$ and the bottom triangle commutes.

    To characterize the composite morphism of the left rectangle, we may equivalently characterize the composite map $\pi_U \sigma_2 h \sigma_1 f\colon \prod\limits_U^S X_U \rightarrow \CoInd_U^V X_U^{\times \Res_U^V S}$;
    in fact, under the expression $X_U^{\times \Res_U^V S} \simeq \prod\limits_W^{\Res_U^V S} \Res_W^U X_U$, it suffices to characterize the composite morphism $\prod_U^S X_U \rightarrow \CoInd_W^V \Res_W^U X_U$ and verify that it is homotopic to the relevant projection of $g$ for each $W,U$.
    
    In particular, the relevant projection of $g$ is the composite morphism
    \[
      \prod_U^S X_U \twoheadrightarrow  \CoInd_U^V X_U \xrightarrow{\delta_{U,W}} \CoInd_W^V \Res_W^U X_U
    \]
    where $\delta_{U,W}$ is a Kronecker delta\footnote{Here, $0$ is the image under $U$ of the zero map in the (pointed) $\infty$-category $\Mon_{\Res_U^V \cO}(\Res_U^V \cC)$.}
    \[
      \delta_{U,W} = \begin{cases}
        \id & U = W; \\ 
        0 & \mathrm{otherwise}. 
      \end{cases}  
    \]
    Moreover, note that the projection $\pi_U \sigma_2 h \sigma_1\colon \prod\limits_U^S X_U \rightarrow X_U^{\times \Res_U^V S}$ itself factors as
    \[
      \prod\limits_U^S \prn{\Res_U^V \coprod_U^V X_U} \twoheadrightarrow \CoInd_U^V X_U \xrightarrow{\;\; \widetilde f_U \;\;} \CoInd_W^V \Res_W^U X_U,
    \]
    so we're tasked with verifying that $\widetilde f_U$ is homotopic to $\delta_{U,W}$.
    Indeed, this follows by examining the following diagram:
  \[
    \begin{tikzcd}[ampersand replacement=\&, row sep=small]
      {\prod\limits_U^S X_U} \& {\prod\limits_U^S \prn{\Res_U^V\coprod\limits_U^V X_U}} \& {\prn{\coprod\limits_U^S X_U}^{\times S}} \& {\prn{\prod\limits_U^S X_U}^{\times S}} \& {\prod\limits_U^S X_U^{\times \Res_U^V S}} \\
      {\CoInd_U^V X_U} \&\& {\CoInd_U^V \Res_U^V\coprod\limits_U^V X_U} \& {\CoInd_U^V \Res_U^V\prod\limits_U^V X_U} \& {X_U^{\times \Res_U^V S}} \\
      \&\&\& {\CoInd_W^V \Res_W^U X_U}
      \arrow["{{f}}", from=1-1, to=1-2]
      \arrow[two heads, from=1-1, to=2-1]
      \arrow["\simeq"{marking, allow upside down}, draw=none, from=1-2, to=1-3]
      \arrow[two heads, from=1-2, to=2-3]
      \arrow["h", from=1-3, to=1-4]
      \arrow["\simeq"{marking, allow upside down}, draw=none, from=1-4, to=1-5]
      \arrow[two heads, from=1-5, to=2-5]
      \arrow["{\CoInd_U^V i_U}", from=2-1, to=2-3]
      \arrow["{\delta_{U,W}}"{description}, curve={height=12pt}, from=2-1, to=3-4]
      \arrow["{\CoInd_U^V W}", from=2-3, to=2-4]
      \arrow["\simeq"{marking, allow upside down}, draw=none, from=2-4, to=2-5]
      \arrow[two heads, from=2-4, to=3-4]
    \end{tikzcd}
  \] 
\end{proof}

\begin{proof}[Proof of \cref{k-connected theorem}]
  Assume $\cO^{\otimes}$ is $\ell$-connected at $I$, i.e. \cref{cond:ell-connected at I}.
  Writing $X$ for an $S$-tuple in $\uMon_{\cO}(\cC)$, we chase $W_{S,X}$ around the following commutative diagram:
  \[
    \begin{tikzcd}
      \uMon_{\cO}(\tau_{\leq \ell} \cC) \arrow[r, "\iota_{\cO}"] \arrow[d, "U_{\leq \ell}"]
      & \uMon_{\cO}(\cC) \arrow[r, "\tau_{\cO}"] \arrow[d, "U"]
      & \uMon_{\cO}(\tau_{\leq \ell}\cC) \arrow[d, "U_{\leq \ell}"]\\
      \uCoFr^{\cT} \tau_{\leq \ell} \cC \arrow[r, "\iota"]
      & \uCoFr^{\cT} \cC \arrow[r, "\tau_{\leq \ell}"]
      & \uCoFr^{\cT} \tau_{\leq \ell} \cC
    \end{tikzcd}
  \]
  In particular, by \cref{hn smashing corollary,Wirthmuller truncation lemma}, $\tau_{\cO} W_{S,X} \sim W_{S, \tau_{\cO} X}$ is an equivalence, so
  \[
    U_{\leq \ell} \tau_{\cO} W_{S,X} \sim \tau_{\leq \ell} U W_{S,X}
  \]
  is an equivalence, i.e. $U W_{S,X}$ is an $\ell$-equivalence.
  In turn, by \cref{section lemma} this implies that $U W_{S,X}$ is $\ell$-connected, so \cref{U conservativity} implies that $W_{S,X}$ is $\ell$-connected, i.e. \cref{cond:topos wirthmuller maps}.
  
  The implication \cref{cond:topos wirthmuller maps} $\implies$ \cref{cond:spaces wirthmuller maps} is immediate, so assume \cref{cond:spaces wirthmuller maps}, i.e. fix the case $\cC \deq \cS$ and assume that that $W_{S,X}$ is $\ell$-connected for all $X \in \Alg_{\cO} \cS$ and $S \in \uFF_I$.
  We may invert the above argument:
  this time, since $\tau_{\cO}$ is a left adjoint, it preserves $\ell$-connectedness, so $\tau_{\cO} W_{S, X}$ is $\ell$-connected.
  Hence
  \[
    U_{\leq \ell} \tau_{\cO} W_{S, X} = \tau_{\leq \ell} U W_{S, X}
  \]
  is an equivalence, i.e. $UW_{S,X}$ is an $\ell$-equivalence.
  We've shown that \cref{cond:spaces wirthmuller maps} implies the following.
  \begin{enumerate}[label={(\alph*)}]\setcounter{enumi}{3}
    \item \label[condition]{cond:spaces wirthmuller maps weak} For all 
      $I$-admissible $V$-sets $S \in \FF_{I,V}$ and
      $S$-tuples of $\cP$-monoids
      $X = (X_K) \in \uMon_{\cP}(\cS)_S$,
      the underlying $V$-space map of the $S$-indexed $\cP$-monoid space Wirthm\"uller map
      \[
        UW_{S, X}\colon U\coprod_K^S X_K \longrightarrow U\prod_K^S X_K \simeq \prod_K^S U X_K 
      \]
      is an $\ell$-equivalence.
  \end{enumerate}
  
  We're left with showing that \cref{cond:spaces wirthmuller maps weak} implies \cref{cond:ell-connected at I};
  indeed, choosing $X = \iota_{\cO} Y$, the map 
  \[
    \tau_{\leq \ell} UW_{S,\iota_{\cO}Y } = U_{\leq \ell} \tau_{\cO} W_{S,\iota_{\cO} Y} = U_{\leq \ell} W_{S, \tau_{\cO} \iota_{\cO} Y} = W_{S,Y}
  \]
  is an equivalence, i.e. $\uMon_{\cO}(\tau_{\leq \ell} \cC)$ is $I$-semiadditive.
  Thus \cref{Mapping spaces prop} implies \cref{cond:ell-connected at I}.
\end{proof}

\subsection{The proof of \texorpdfstring{\cref{LEHA,HEHA}}{the main theorems}}
\label{Reduction subsection}
Note the following.
\begin{proposition}\label{Reduced endomorphism operad corollary}
  If $\cP^{\otimes}$ is $\ell$-connected at $I$, then for all $(k + \ell + 2)$-toposes $\cC$, the reduced endomorphism $I$-operad $\End_X\prn{\uMon_{\cP}(\cC)^{I-\times}}$ is an $I$-$(k + 1)$-operad.
\end{proposition}
\begin{proof}
  Since $\cC$ is a $(k + \ell + 2)$-category, $X$ is $(k + \ell + 2)$-truncated, and \cref{k-connected theorem} implies that $W_{X,S}$ is $\ell$-connected, so the result follows from \cref{Schlank connectivity prop}.
\end{proof}
We quickly acquire a slightly weakened version of \cref{LEHA}. 
\begin{corollary}\label{LEHA easy}\label{HEHA easy}
  Suppose $\cT$ is an atomic orbital $\infty$-category, $\cO^{\otimes}$ and $\cP^{\otimes}$ are unital $\cT$-operads and $I$ is a unital weak indexing category.
If $\cO^{\otimes}$ is $k$-connected at $I$ and $\cP^{\otimes}$ is $\ell$-connected at $I$, then $\cO^{\otimes} \obv \cP^{\otimes}$ is $(k + \ell + 2)$-connected at $I$.
\end{corollary}
  \begin{proof}
  By \cref{Reduced endomorphism operad corollary}, we know that 
  $\End^{I,\red}_X\prn{\uMon_{\cP}(\cC)^{\cT-\times}} \simeq \End^{I,\red}_X \prn{\uMon_{\cP}(\cC)^{I-\times}}$ is an $I$-$(k+1)$-operad for $\cC$ a $(k + \ell + 2)$-topos;
  by \cref{hn smashing corollary} this shows that $\Mon_{\cO} \End_X \prn{\uMon_{\cP}(\cC)^{\cT-\times}}$ is $I$-cocartesian, so \cref{Checking morita equivalences on End} shows that
  \[
    \CMon_{I} \Mon_{\cO} \Mon_{\cP}(\cC)^{\simeq} \xrightarrow{\;\; \sim \;\;} \Mon_{\cO} \Mon_{\cP}(\cC)^{\simeq}.
  \]
  By \cref{Mapping spaces prop}, this implies that the map 
  \[
    \cO^{\otimes} \obv \cP^{\otimes} \simeq \cO^{\otimes} \obv \cP^{\otimes} \obv \triv_{\cT}^{\otimes} \xrightarrow{\;\;\;\;\; \id \otimes \id \otimes ! \;\;\;\;\;} \cO^{\otimes} \obv \cP^{\otimes} \obv \cN_{I\infty}^{\otimes}
  \]
  is an $h_{k + \ell + 2}$-equivalence, so \cref{hn smashing corollary} shows that $\cO^{\otimes} \obv \cP^{\otimes}$ is $\prn{k + \ell + 2}$-connected at $I$.
\end{proof}
From this, we conclude the main theorems of this paper.
\begin{proof}[Proof of \cref{HEHA,LEHA}]
  Restriction assembles to a (tautologically symmetric monoidal) equivalence
  \[
    \Op_\cT^{\otimes} \simeq \lim_{V \in \cT} \Op_V^{\otimes}
  \]
  such that, given a morphism $V \rightarrow W$ in $\cT$ and $S$ a finite $V$-set, $\Res_V^{\cT} \cO(S) \simeq \cO(S)$.
  In particular, \cref{LEHA,HEHA} may be verified after restriction to each to $V \in \cT$, in which case the base $\infty$-category $\cT_{/V}$ has a terminal object.
  
  Moreover, each $\cO(S)$ and $\cP(S)$ are easily determined by arity support except in the case $V \in \upsilon(\cO) = \upsilon(\cP)$, and arity support is additive in the predicted way by \steprop{1.44};
  thus \cref{LEHA,HEHA} may be verified after restriction to each $V \in \upsilon(\cO)$, in which case $\Res_V^{\cT} \cO^{\otimes}$ and $\Res_V^{\cT} \cP^{\otimes}$ are unital.
  This and \cref{LEHA easy} together yield \cref{LEHA}, and \cref{HEHA} follows by setting $I \deq A\cO$.
\end{proof}

\subsection{Sharpness}  
In this subsection, we show that the inequalities of \cref{HEHA,LEHA} are not always attained.
One reason for this is the discrepancy between unions and joins of weak indexing systems.
\begin{example}
      It follows by definition that
  \[
    \Conn_{\cN_{I \infty}}(J) = \begin{cases}
      \infty & J \subset I, \text{ and} \\ 
      -2 & \mathrm{otherwise}; 
    \end{cases}
  \]
  we also found in \cite{Tensor} that $\cN_{I \infty}^{\otimes} \obv \cN_{J \infty}^{\otimes} \simeq \cN_{I \vee J\infty}^{\otimes}$.
  Generically, this defeats sharpness of \cref{LEHA}, as
  \[
    \prn{\Conn_{\cN_{I \infty}} + \Conn_{\cN_{J \infty}} + 2}^{-1}(\infty) = \wIndSys^{a\uni}_{\cT, \leq I} \cup \wIndSys^{a\uni}_{\cT, \leq J} \subsetneq \wIndSys^{a\uni}_{\cT, \leq I \vee J} = \Conn_{\cN_{I\infty} \otimes \cN_{J \infty}}^{-1}(\infty).\qedhere
  \]
\end{example}
Another reason is topological;
in forthcoming work, given $V$ an orthogonal $G$-representation, we will show that the little $V$-disks $G$-operad $\EE^{\otimes}_V$ is $\ell$-connected at the minimal unital weak indexing category $I_S \vee I^0$ containing $S$ if and only if the following conditions are satisfied:
\begin{enumerate}[label={(\alph*)}]
  \item For all orbits $[H/K] \subset S$ and intermediate inclusions $K \subset J \subset H$, we have $\dim V^J \geq \dim V^K + \ell + 2$, and
  \item if $\abs{S^H} \geq 2$, then $\dim V^H \geq \ell + 2$.
\end{enumerate}
Moreover, we will show that $\EE_V$ is additive under tensor products, i.e. $\EE^{\otimes}_{V} \obv \EE^{\otimes}_W \simeq \EE^{\otimes}_{V \oplus W}$.
\begin{example}
  Let $G \deq C_2$, with sign representation $\sigma$.
  Then, we have fixed point dimensions
  \[
    \dim \prn{a + b\sigma}^{e} = a + b; \hspace{40pt} \dim \prn{a + b\sigma}^{c_2} = a.
  \]
  In particular, the connectivity function has
  \begin{align*}
    \Conn_{\EE_{a + b\sigma}}(k *_e) &= a + b - 2\\
    \Conn_{\EE_{a + b\sigma}}(c *_{C_2} + d[C_2/e])) &= \begin{cases}
      a-2 & d = 0 \\ 
      b-2 & c < 2 \\
      \min(a,b) - 2 & \mathrm{otherwise}.
    \end{cases}
  \end{align*}
  where $\Conn(S) \deq \Conn(I_S \vee I^0)$.
  Note that $\Conn_{\EE_{a + b\sigma}}(c *_{C_2} + d[C_2/e])$ is as non-additive as is possible in the last case;
  indeed, the examples $1 + b\sigma$ and $a' + \sigma$ have the same arity-support, but when $a',b > 1$, we have
  \begin{align*}
    \Conn_{1 + b\sigma}(2*_{C_2} + [C_2/e]) + 
    \Conn_{a' + \sigma}(2*_{C_2} + [C_2/e]) - 2 &= 0\\
    &< \min(a',b) - 1\\
    &= \Conn_{a' + 1 + (b + 1)\sigma}(2*_{C_2} + [C_2/e]).\qedhere
  \end{align*}
\end{example}

\section{Corollaries}
\subsection{Smashing localizations and \texorpdfstring{\cref{cor:smashing}}{Corollary E}} 
\cref{HEHA} specializes to infinite tensor powers as follows.
\begin{corollary}\label{Classification of idempotents}
  Suppose $\cO^{\otimes}$ is an almost-reduced $\cT$-operad.
  Then, the following conditions are equivalent.
  \begin{enumerate}[label={(\alph*)}]
    \item \label[condition]{cond:auw} $\cO^{\otimes}$ is an almost-unital weak $\cN_\infty$-operad.
    \item \label[condition]{cond:EHA} ($\cO^{\otimes}$-EHA) the unique map $\triv_{\cT}^{\otimes} \rightarrow \cO^{\otimes}$ yields an equivalence
      \[
        \cO^{\otimes} \simeq \cO^{\otimes} \obv \triv_{\cT}^{\otimes} \xrightarrow{\; \id \otimes ! \;} \cO^{\otimes} \obv \cO^{\otimes}.
      \]
    \item \label[condition]{cond:AEHA} (abstract $\otimes$-idempotence) there exists an equivalence $\cO^{\otimes} \obv \cO^{\otimes} \simeq \cO^{\otimes}$.
  \end{enumerate}
\end{corollary}
\begin{proof}
  The implication \ref{cond:auw} $\implies$ \ref{cond:EHA} is \stecor{2.4}, and is also implied by \cref{HEHA}.
  The implication \ref{cond:EHA} $\implies$ \ref{cond:AEHA} is obvious.
  To see the implication \ref{cond:AEHA} $\implies$ \ref{cond:auw}, note that \cref{HEHA} implies that $\cO^{\otimes}$ is $\infty$-connected, i.e. all of its nonempty structure spaces are contractible.
  The result follows by the identification of such almost-reduced $\cT$-operads with almost-unital weak $\cN_\infty$-operads 
  \cite[Thm.~\href{https://arxiv.org/pdf/2501.02129v1\#cooltheorem.3}{C}]{EBV}.
\end{proof}

\begin{remark}
To see why \cref{cond:EHA} is an \emph{Eckmann-Hilton argument}, note that it is equivalent to the condition that $\cO^{\otimes}$ possesses a unital magma structure in $\Op_{\cT}^{\otimes}$ whose multiplication map $\mu\colon \cO^{\otimes} \obv \cO^{\otimes} \rightarrow \cO^{\otimes}$ is an equivalence;
unitality of $\mu$ is precisely the condition that the pullback natural transformation
\[
  \delta = \mu^*\colon \Alg_{\cO}(\cC) \longrightarrow \Alg_{\cO} \uAlg_{\cO}^{\otimes}(\cC)
\]
is split by restriction to either $\cO$-algebra structure, and the fact that $\mu$ is an equivalence is precisely the condition that $\delta$ is a natural equivalence, i.e. pairs of interchanging $\cO$-algebra structures agree and $\cO$-algebra structures interchange with themselves in an essentially unique way.
\end{remark}

Now, we may use this to prove \cref{cor:smashing}.
\begin{proof}[Proof of \cref{cor:smashing}]
  We proved in 
  \cite[Thm.~\href{https://arxiv.org/pdf/2501.02129v1\#cooltheorem.3}{C}]{EBV}
  that weak $\cN_\infty$-operads are equivalent to subterminal $\cT$-operads, yielding the horizontal equivalences.
  
  Now, we showed in \sterem{2.33} that every idempotent algebra in $\Op_G^{a\uni}$ is almost-reduced, so for the vertical equivalences, we may replace $\Op_G^{a\uni}$ with $\Op_G^{a\red}$;
  in this symmetric monoidal category, the unit $\EE_0^{\otimes}$ is initial, so being an idempotent algebra is a property.
  Indeed, \cref{Classification of idempotents} shows that this is the same property as being subterminal.
  Thus the vertical equivalences follow from the correspondence between idempotent algebras and smashing localizations 
  (see \cite[\S~\href{https://arxiv.org/pdf/1305.4550v1\#section.3}{3}]{Gepner} and \cite[\S~\href{https://arxiv.org/pdf/2007.13089v2\#subsection.5.1}{5.1}]{Carmeli}).
  
  What's left is identifying the diagonal arrow, i.e. characterizing the smashing localization associated with $\cN_{I\infty}^{\otimes}$;
  this was done in \stethm{2.6}.
\end{proof}

\subsection{Equivariant loop spaces and \texorpdfstring{\cref{I-spectrum coolcorollary}}{Cor F}}
Let $G$ be a finite group and $I$ an indexing category.
\begin{lemma}\label{SpI algebraic}
  If $k \geq 0$, then there exists an equivalence $\Sp_{I, [k,\ell]} \simeq \CAlg_{I}(\cS_{G, [k,\ell]})$ over $\cS_{G, [k,\ell]}$.
\end{lemma}
\begin{proof}
  Given a model $\cO^t \in \Op(\sSet_G)$ for $\cN_{A\cO}^{\otimes}$, we acquire a diagram of equivalences
  \[\begin{tikzcd}[ampersand replacement=\&, column sep=small]
	{\Sp_{A\cO, [k,\ell]}} \& {\Alg_{\cO^t}\prn{\Top_{G, [k,\ell]}}\brk{\mathrm{WEQ}^{-1}}} \&\&\& {\CMon_{A\cO}\prn{\cS_{[k,\ell]}}} \&\&\& {\CAlg_{A\cO}\prn{\ucS_{G, [k,\ell]}}} \\
	\& {\cS_{G, [k,\ell]}}
	\arrow["\deq"{description}, draw=none, from=1-1, to=1-2]
	\arrow["{\Omega^\infty}" below, from=1-1, to=2-2]
	\arrow["{\text{\cite{Marc}}}", "\sim" below, from=1-2, to=1-5]
	\arrow["U"', from=1-2, to=2-2]
	\arrow["{\text{\cite{Tensor}}}", "\sim" below, from=1-5, to=1-8]
	\arrow["U"{description}, from=1-5, to=2-2]
	\arrow["U", from=1-8, to=2-2]
\end{tikzcd}\]
Indeed, the case without connectivity and truncations
is proved drectly in the cited articles, and the restriction to $[k,\ell]$ follows by unwinding definitions to see
that the notion of concentration at degrees $[k,\ell]$ corresponds with the preimage of $\cS_{G, [k,\ell]} \subset \cS_G$ within each category.
\end{proof}
From this, we're ready to conclude \cref{I-spectrum coolcorollary}.
\begin{proof}[Proof of \cref{I-spectrum coolcorollary}]
  First note that the $G$-$\infty$-category of $k$-connected $\ell$-truncated connected $G$-spaces is a $G$-$(\ell - k)$-category;
  indeed, Elmendorf's theorem yields an equivalence
  \[
    \prn{\ucS_{G, [k,\ell]}}_H \simeq \Fun\prn{\cO_H^{\op}, \cS_{[k,\ell]}},
  \]
  and $\cS_{[k,\ell]}$ is an $(\ell - k+1)$-category as whenever $X$ is $(k-1)$-connected and $Y$ is $\ell$-truncated, we have 
  \[
    \Omega^{\ell - k + 1} \Map(X,Y) \simeq \Map(\Sigma^{\ell - k + 1} X, Y) \simeq *;
  \]
  hence \httcor{2.3.4.8} implies that each value $\prn{\ucS_{G,[k,\ell]}}_H$ is an $(\ell - k+1)$-category.
  Thus \cref{Connectivity corollary,SpI algebraic}, 
  together construct the desired equivalence
  \[\begin{tikzcd}[ampersand replacement=\&, column sep=tiny]
    {\Sp_{A\cO, [k,\ell]}} \& {\CAlg_{A\cO}\prn{\ucS_{G, [k,\ell]}}} \&\&\& {\overbrace{\Alg_{\cO} \uAlg^{\otimes}_{\cO} \cdots \uAlg^{\otimes}_{\cO}}^{(\ell - k + 2)\text{-fold}}\prn{\ucS_{G, [k,\ell]}},} \\
    \& {\cS_{G, [k,\ell]}}
    \arrow["\simeq"{marking, allow upside down}, draw=none, from=1-1, to=1-2]
    \arrow["{\Omega^\infty}"', from=1-1, to=2-2]
    \arrow["U", from=1-2, to=1-5, "\sim" below]
    \arrow["U"{description}, from=1-2, to=2-2]
    \arrow["U", from=1-5, to=2-2]
  \end{tikzcd}\]
\end{proof}

\begin{remark}
  We chose to specialize to the connected setting for convenience;
  one could instead assume that there exists some $\mu \in \cO(2 \cdot *_G)$ whose action on one of the $\cO$-structures on $X$ induces an \emph{invertible} magma structure on the coefficient system $\upi_0 X$, in which case the corresponding $A\cO$-commutative algebra has an underlying grouplike commutative monoid structure; 
  the variation of \cref{I-spectrum coolcorollary} follows \emph{mutatis mutandis}.
\end{remark}

More traditionally, we acquire $\Omega^V$-spectrum structures in a similar circumstance.
\begin{corollary}
  Fix $V$ an orthogonal $G$-representation, $0 \leq k \leq \ell \leq \infty$ related numbers, and $\cO^{\otimes}$ an almost-reduced $G$-operad with $\cO(S) \neq \emptyset$ whenever there exists an embedding $S \hookrightarrow \Res_H^G V$.
  If $X$ is a $(k-1)$-connected and $\ell$-truncated $G$-space admitting $(\ell - k + 2)$-many interchanging $\cO$-algebra structures,
  then $X$ underlies a $V$-infinite loop space.
\end{corollary}
\begin{proof}
  The desired $V$-infinite loop space structure corresponds under the recognition principle of \cite{Rourke,Guillou-May,Juran} with the $\EE_{\infty V}$-structure pulled back along the unique map specified by \cref{Connectivity corollary}:
  \[
    \EE^{\otimes}_{\infty V} \simeq \cN_{AV}^{\otimes} \xrightarrow{\;\;\;\; ! \;\;\;\;} \cN_{A\cO}^{\otimes} \simeq h_{\ell - k + 1}\cO^{\otimes (\ell - k + 2)}.\qedhere
  \]
\end{proof}

\subsection{The \texorpdfstring{$C_p$-operad $\AA^{\otimes}_{2,C_p} \obv \AA^{\otimes}_{2,C_p}$}{Cp-operad A2-Cp tensor A2-Cp} and \texorpdfstring{\cref{GEHA}}{Theorem A}}\label{A2 section}
For the rest of this article, we specialize to $G = C_p$ the group of prime order $p$ and $\cC$ a 1-category.
As in \cref{thm:free operads}, let $\Fr_{\Sigma}(S)$ denote the free $C_p$-symmetric sequence on an operation in arity $S$.
Now, the pointwise formula for left Kan extensions yields equivalences
\begin{equation}\label{eq:Sigma p}
  \begin{split}
  \Fr_{\Sigma,p \cdot *_{C_p}}(*)(p \cdot *_e) &\simeq \colim_{\Res_e^{C_p} p \cdot *_{C_p} \xrightarrow{\sim} p \cdot *_e} * \simeq \Sigma_p;\\
  \Fr_{\Sigma, [C_p/e]}(*)(p \cdot *_e) &\simeq \colim_{\Res_e^{C_p} [C_p / e] \xrightarrow{\sim} p \cdot *_e} * \simeq \Sigma_p.
\end{split}
\end{equation}
We define the $C_p$-symmetric sequence of sets $F_{2,C_p}$ as the coequalizer
\[
  F_{2,C_p} \deq \mathrm{CoEq}
  \prn{\Sigma_p[p \cdot *_e] \rightrightarrows \prn{\Fr_{\Sigma,[C_p/e]}(*) \sqcup \Fr_{\Sigma, p \cdot *_{C_p}(*)}}},
\]
where $\Sigma_p [p \cdot *_e]$ is the $C_p$-symmetric sequence defined by
\[
  \Sigma_p [p \cdot *_e](S) \deq \begin{cases}
    \Sigma_p & S = p \cdot *_e; \\ 
    \emptyset & \mathrm{otherwise}. 
  \end{cases}
\]
and the two arrows are the inclusions of $\Sigma_p[p \cdot *_e]$ into the summands prescribed by \cref{eq:Sigma p}.
We define the unital $C_p$-operad $\AA_{2,C_p}^{\otimes}$ by the Boardman-Vogt tensor product
\[
  \AA_{2,C_p}^{\otimes} \deq \EE^{\otimes}_0 \obv \Fr_{\Op} \prn{F_{2,C_p}}.
\]
As promised, we verify that $\AA_{2,C_p}$-monoids are the same as $C_p$-unital magmas.
\begin{proposition}\label{thm:cp unital magmas}
  There is an equivalence between $\Mon_{\AA_{2,C_p}}(\cC)$ and $C_p$-unital magmas in $\cC$.
\end{proposition}
\def\TX{\widetilde{X}}
\begin{proof}
  By \cref{ex:e0,thm:unital prop} we have
  \[
    \Mon_{\AA_{2,C_p}}(\cC) \simeq \Mon_{\Fr_{\Op}(F_{2,C_p})} \uMon_{\EE_0}^{\otimes}(\cC) \simeq \Mon_{\Fr_{\Op}(F_{2,C_p})} \cC_*.
  \]
  Moreover, by \cref{thm:1-categorical algebras}, the data of an $\AA_{2,C_p}$-monoid structure on $X \in \CoFr^{C_p} \cC$ is equivalently viewed as a map $\eta\colon *_{C_p} \rightarrow X$ (which we identify with an element $\widetilde X \in \CoFr^{C_p}\cC_{*}$) and an element of
  \begin{align*}
    \Mon_{\Fr_{\Op}(F_{2,C_p})}(\End_{\TX}(\cC_*))^{\simeq} 
    &\simeq \Hom_{\Fun\prn{\Tot \uSigma_{C_p}, \cS}}\prn{F_{2,C_p}, \End_{\TX}(\cC_*)}\\
    &\simeq \Hom_{\CoFr^{C_p} \cC_*}\prn{\TX^p, \TX} \times_{\Hom_{\cC_*}\prn{\prn{\TX^e}^p, \TX^e}} \Hom_{\CoFr^{C_p} \cC_*}\prn{\CoInd_e^{C_p} \TX^e, \TX}.
  \end{align*}
  We're left with interpreting this concretely:
  by a standard argument, $\Hom_{\CoFr^{C_p} \cC_*}(\TX^p, \TX)$ corresponds bijectively with the set of unital magma structures on $X$ with unit $\eta$, and this corresponds bijectively with the pairs of unital magma structures on $X^{C_p}$ and $X^e$ with unit maps $\eta^{C_p}$ and $\eta^e$ such that the restriction map is a homomorphism.
  Under this bijection, the forgetful map $\Hom_{\CoFr^{C_p} \cC_*} \prn{\TX^p, \TX} \rightarrow \Hom_{\cC_*} \prn{\prn{\TX^e}^p, \TX}$ simply forgets the data of $X^{C_p}$ and the restriction.

  Similarly, since $C_p$-coefficient coinduction is presented by the coefficient system $X^p \xleftarrow{\Delta} X$ with permutation action, $\Hom_{\CoFr^{C_p} \cC^*}\prn{\CoInd_e^{C_p} \TX^e, \TX}$ corresponds bijectively with the set of unital $C_p$-equivariant transfers $t\colon X^{e} \rightarrow X^{C_p}$ and unital magma structures on $X^e$ with unit $\eta^e$ satisfying the condition that the following diagram commutes.
  \[
     \begin{tikzcd}[ampersand replacement=\&]
      {X^e} \& {X^{C_p}} \\
      {\prn{X^e}^p} \& {X^e}
      \arrow["t", from=1-1, to=1-2]
      \arrow["\Delta", from=1-1, to=2-1]
      \arrow["r", from=1-2, to=2-2]
      \arrow["*", from=2-1, to=2-2]
    \end{tikzcd}
  \]
  Once again, the forgetful map restricts to the unital magama structure on $\eta^e$;
  thus the fiber product corresponds exactly with $G$-unital magma structures on $X$ with units $\eta^e$ and $\eta^{C_p}$.

  Now, what we've described is a bijective assignment of \emph{sets} $\Ob \Mon_{\AA_{2,C_p}}(\cC) \rightarrow \Ob \mathrm{Magma}^{\uni}_{C_p}(\cC)$ over $\Ob \cC$.
  To conclude, it suffices to prove that a $\CoFr^{C_p}\cC$ morphism between a pair of $C_p$-unital magmas is a $C_p$-unital magma homomorphism if and only if it's an $\AA_{2,C_p}$-algebra homomorphism.

  To prove this, note that an $\AA_{2,C_p}$-monoid morphism is equivalently a $\Fr_{\Op}(F_{2,C_p})$-monoid morphism of pointed objects, i.e. a pair of maps $F^e\colon M^e \rightarrow N^{e}$ and $F^{C_p}\colon M^{C_p} \rightarrow N^{C_p}$ which are compatible with units, satisfying $F^{C_p} \circ t = t \circ F^e$ and $F^e \circ r = r \circ F^{C_p}$ together with $p$-degree additivity 
  \[\begin{tikzcd}[ampersand replacement=\&]
	{\prn{M^{C_p}}^p} \& {\prn{N^{C_p}}^p} \& {\prn{M^e}^p} \& {\prn{M^e}^p} \\
	{M^{C_p}} \& {N^{C_p}} \& {M^e} \& {N^e}
	\arrow[from=1-1, to=1-2]
	\arrow[from=1-1, to=2-1]
	\arrow[from=1-2, to=2-2]
	\arrow[from=1-3, to=1-4]
	\arrow[from=1-3, to=2-3]
	\arrow[from=1-4, to=2-4]
	\arrow[from=2-1, to=2-2]
	\arrow[from=2-3, to=2-4]
\end{tikzcd}\]
  It suffices to note that a map between the pointed sets underlying unital magmas is a homomorphism if and only if it intertwines with $n$ary addition for \emph{some} $n \geq 2$;
  indeed, one can simply identify binary addition with $n$-ary addition whose first $(n-2)$-factors are the unit.
\end{proof}

We now spell out the interchange relations explicitly.
\begin{proposition}
  There is an equivalence between $\Mon_{\AA_{2,C_p} \otimes \AA_{2,C_p}}(\cC)$ and pairs of $G$-unital magma structures $(M,*,\bullet,t_*,t_{\bullet})$ in $\cC$ satisfying the interchange relations
  $1_* = 1_\bullet$ and
  \[\begin{tikzcd}[ampersand replacement=\&, column sep=small]
  {\prn{X^{p}}^p} \& {X^p} \& {X^{C_p}} \& {X^e} \& {X^{C_p}} \& {\prn{X^e}^p} \& {\prn{X^{C_p}}^p} \& {\prn{X^e}^p} \& {\prn{X^{C_p}}^p} \\
	{X^p} \& X \& {X^e} \& {\prn{X^e}^p} \& {X^e} \& {X^e} \& {X^{C_p}} \& {X^e} \& {X^{C_p}}
  \arrow["\prn{\bullet}", from=1-1, to=1-2]
  \arrow["{\prn{*}}"', from=1-1, to=2-1]
	\arrow["{*}"', from=1-2, to=2-2]
	\arrow["r"', from=1-3, to=2-3]
	\arrow["{t_\bullet}"', from=1-4, to=1-3]
	\arrow["{t_*}", from=1-4, to=1-5]
	\arrow["\Delta"', from=1-4, to=2-4]
	\arrow["r", from=1-5, to=2-5]
	\arrow["{\prn{t_\bullet}}", from=1-6, to=1-7]
	\arrow["{*}"', from=1-6, to=2-6]
	\arrow["{*}"', from=1-7, to=2-7]
	\arrow["{\prn{t_*}}", from=1-8, to=1-9]
	\arrow["\bullet"', from=1-8, to=2-8]
	\arrow["\bullet"', from=1-9, to=2-9]
	\arrow["\bullet", from=2-1, to=2-2]
	\arrow["{*}", from=2-4, to=2-3]
	\arrow["\bullet"', from=2-4, to=2-5]
	\arrow["{t_\bullet}"', from=2-6, to=2-7]
	\arrow["{t_*}"', from=2-8, to=2-9]
\end{tikzcd}\]
\end{proposition}
\begin{proof}
  \cref{ex:e0,thm:unital prop} yields an equivalence.
  \[
    \Mon_{\AA_{2, C_p}^{\otimes 2}}(\cC) \simeq \Mon_{\Fr_{\Op}(F_{2,C_p})^{\otimes 2}}(\cC_*).
  \]
  This is characterized explicitly by \cref{thm:1-categorical interchange,thm:cp unital magmas};
  it suffices to note that the specified interchange relations correspond precisely with the conditions that $t_{\bullet}$ and $\bullet$ are $C_p$-unital magma homomorphisms.
\end{proof}
 
We conclude the following form of \cref{GEHA}.
\begin{corollary}\label{GEHA-precise}
  Given $\cC$ a 1-category, the forgetful functor
  \begin{align*}
  \Fun^{\times}(\Span(\FF_{C_p}), \cC) \longrightarrow& \Mon_{\AA_{2,C_p} \otimes \AA_{2,C_p}}(\cC)\\
    \simeq& \cbr{\text{Interchanging pairs of } C_p \text{-unital magmas in } \cC}
  \end{align*}
  is an equivalence of categories.
\end{corollary}

\printbibliography

\end{document}